\documentclass[12pt]{amsart}
\usepackage{hero}
\usepackage{turnstile}
\usepackage{psfrag,epsfig,multimedia}
\usepackage{graphicx}
\usepackage{amssymb,amsmath,amsfonts}

\usepackage[tmargin=1.15in, bmargin=1.15in, rmargin=1.0in, lmargin=1.0in]{geometry}

\def\ra{\mbox{\rm a}}
\def\etal{{\it {et al}}}

\def\ra{\mbox{\rm a}}
\def\mX{{\mathbb X}}

\def\Beta{B}

\begin{document}


\title{Hub discovery in partial correlation graphical models\protect}


\author{Alfred Hero}
\address{Departments of EECS, BME and Statistics,
University of Michigan - Ann Arbor, MI 48109-2122, U.S.A}
\email{hero@umich.edu}

\author{Bala Rajaratnam}
\address{Department of Statistics,
Stanford, CA 94305-4065, U.S.A.}
\email{brajarat@stanford.edu}

\maketitle

\begin{abstract}


One of the most important problems in large scale inference problems
is the identification of variables that are highly dependent on
several other variables. When dependency is measured by partial
correlations these variables identify those rows of the partial
correlation matrix that have several entries with magnitude close to
one; i.e., hubs in the associated partial correlation graph.  This
paper develops theory and algorithms for discovering such hubs from a few observations of these variables.  We introduce a hub
screening framework in which the user specifies both a minimum
(partial) correlation $\rho$ and a minimum degree $\delta$ to screen
the vertices. The choice of $\rho$ and $\delta$ can be guided by our
mathematical expressions for the phase transition correlation
threshold $\rho_c$ governing the average number of discoveries. We
also give asymptotic expressions for familywise discovery rates under
the assumption of large $p$, fixed number $n$ of multivariate samples,
and weak dependence. Under the null hypothesis that the covariance
matrix is sparse 
these
limiting expressions can be used to enforce FWER constraints and to
rank these discoveries in order of statistical significance (p-value).
For $n\ll p$ the
computational complexity of the proposed partial correlation
screening method is low and is therefore highly scalable.
Thus it can be applied to significantly
larger problems than previous approaches.
The theory is applied to discovering hubs in a
high dimensional gene microarray dataset.
\end{abstract}


\bc
\bf Keywords
\ec

Gaussian graphical models, correlation networks,
nearest neighbor dependency, node degree and connectivity,  asymptotic Poisson limits, discovery rate phase transitions,  p-value trajectories


\section{Introduction}


This paper treats the problem of screening a $p$-variate sample for
strongly and multiply connected vertices in the partial correlation
graph associated with the the partial correlation matrix of the
sample. This problem, called hub screening, is important in many
applications ranging from network security to computational biology to
finance to social networks. In the area of  network security, a
node that becomes a hub of high correlation with neighboring nodes
might signal anomalous activity such as a coordinated flooding attack.
In the area of computational biology the set of hubs of a gene
expression correlation graph can serve as potential targets for drug
treatment to block a pathway or modulate host response. In the area of
finance a hub might indicate a vulnerable financial instrument or
sector whose collapse might have major repercussions on the market. In
the area of social networks a hub of observed interactions between
criminal suspects could be an influential ringleader. The techniques
and theory presented in this paper permit scalable and reliable
screening for such hubs. Unlike the correlations screening problem studied in \cite{Hero&Rajaratnam:Paper1_Techreport} this paper considers the more challenging problem of partial correlation screening for variables with a high degree of connectivity. In particular we consider 1) extension
to the more difficult problem of screening for partial correlations exceeding a specified
magnitude; 2) extension to screening variables whose vertex degree in
the associated partial correlation  graph\footnote{Partial correlation graphs
are often called concentration graphs.} exceeds a specified
degree.

The hub screening theory presented in this paper can be applied to
structure discovery in $p$-dimensional Gaussian graphical models, a
topic of recent interest to statisticians, computer scientists
and engineers working in areas such as
gene expression analysis, information theoretic imaging and sensor networks
\cite{friedman2007sic}\cite{liu2010feedback},\cite{wiesel2010covariance}. For example, the
authors of \cite{anandkumar2009detection} propose a Euclidean nearest
neighbor graph method for testing independence in a Gauss-Markov
random field model. When specialized to the null hypothesis of
spatially independent Gaussian measurements, our results characterize
the large $p$ phase transitions and specify a weak (Poisson-like)
limit law on the number of highly connected nodes in such nearest
neighbor graphs for finite number $n$ of observations.

Many different methods for inferring properties of correlation and partial
correlation matrices have been recently proposed
\cite{friedman2008sparse},\cite{rajaratnam2008flexible},
\cite{Rothman&etal:08}, \cite{Bickel&Levina:AOS08_2},
\cite{Peng&etal:JASA09}.  Several of these methods have been
contrasted and compared in bioinformatics applications
\cite{gill2010statistical}, \cite{Kramer&etal:BMC09},
\cite{pihur2008reconstruction} similar to the ones we consider in
Sec. \ref{sec:experiments}.  The above papers address the {\it
  covariance selection} problem \cite{Dempster:Biometrics72}: to find
the non-zero entries in the covariance or inverse covariance matrix
or, equivalently, to find the edges in the associated correlation or
partial correlation graph.


The problem treated and the solution proposed in this paper differ
from those of these previous papers in several important ways: 1) as
contrasted to \cite{Dempster:Biometrics72} our objective is to
screen for vertices in the graph instead of to screen for edges; 2)
unlike 
\cite{Bickel&Levina:AOS08_2} our objective is to directly control false
positives instead of maximizing a likelihood function or minimizing a
matrix norm approximation error; 3) our framework is specifically
adapted to the case of a finite number of samples and a large number
of variables $(n \ll p$); 4) our asymptotic theory provides
mathematical expressions for the p-value for each of the variables
with respect to a sparse null hypothesis on the covariance;
5) unlike lasso type methods like
\cite{Peng&etal:JASA09} the hub screening implementation can be
directly applied to very large numbers of variables without the
need for initial variable reduction. 
Additional relevant literature on correlation based methods can be found in \cite{Hero&Rajaratnam:Paper1_Techreport}.
%

For specified $\rho \in [0,1]$ and $\delta \in 1, \ldots, p-1$, a hub is
defined broadly as any variable that is correlated with at least $\delta$
other variables with magnitude correlation exceeding $\rho$.  Hub
screening is accomplished by thresholding the sample correlation
matrix or partial correlation matrix and searching for rows with more
than $\delta$ non-zero entries.  We call the former {\it correlation hub
  screening} and the latter {\it partial correlation hub screening}. The
screening is performed in a computationally efficient manner by
exploiting the equivalence between correlation graphs and ball graphs
on the set of Z-scores.  Specifically, assume that $n$ samples of $p$
variables are available in the form of a data matrix where $n
<p$. First the columns of the data matrix are converted to standard
$n$-variate Z-scores.  The set of $p$ Z-scores uniquely determine the
sample correlation matrix. If partial correlations are of interest,
these Z-score are replaced by equivalent modified Z-scores that
characterize the sample partial correlation matrix, defined as the
Moore-Penrose pseudo-inverse of the sample correlation matrix.  Then
an approximate k-nearest neighbor (ANN) algorithm is applied to the
Z-scores or the modified Z-scores to construct a ball graph associated
with the given threshold $\rho$.  Hub variables are discovered by
scanning the graph for those whose vertex degree exceeds $\delta$. The ANN
approach only computes a small number of the sample correlations or
partial correlations, circumventing the difficult (or impossible) task
of computing all entries of the correlation matrix when there are
millions (or billions) of variables. State-of-the-art ANN software has
been demonstrated on over a billion variables \cite{Jegou2010product}
and thus our proposed hub screening procedure has potential
application to problems of massive scale. We also note that using
the standard Moore-Penrose inverse is well understood to be a sub-optimal
estimator of the partial correlation matrix in terms of minimum mean square error \cite{GoldsteinSmith:74}. To our knowledge its properties for screening for partial correlations has yet to be investigated. This paper demonstrates that its potential use in detecting high partial correlations as compared to estimating them.

No screening procedure would be complete without error control.  We
establish limiting expressions for mean hub discovery rates.
These expressions are
used to obtain an approximate phase transition threshold $\rho_c$
below which the average number of hub discoveries abruptly
increases. When the screening threshold $\rho$ is below $\rho_c$ the
discoveries are likely to be dominated by false positives. We then
show that the number of discoveries is dominated by a random variable
that converges to a Poisson limit as $\rho$ approaches $1$ and $p$
goes to infinity.  Thus the probability of making at least one hub
discovery converges to a Poisson cumulative distribution.  In the case
of independent identically distributed (i.i.d.) elliptically
distributed samples and sparse block diagonal covariance matrix, the
mean of the dominating variable does not depend on the unknown
correlations.  In this case we can specify asymptotic p-values on
hub discoveries of given degree under a
sparse-covariance null model. Finite $p$ bounds on the  Poisson p-value
approximation error are given that decrease at rates determined by
$p$, $\delta$, $\rho$, and the sparsity factor of the covariance matrix.

To illustrate the power of the hub screening methods we apply them to
a public gene expression datases: the NKI breast cancer data
\cite{chang2005robustness}.
Each of these datasets contains over twenty thousand
variables (genes) but many fewer observations (GeneChips).  The
screening reveals interesting and previously unreported dependency
structure among the variables. For the purposes of exploring
neighborhood structure of the discoveries we introduce a waterfall
plot of their approximate p-values that
plots the family of degree-indexed p-value curves over the range
of correlation thresholds. This graphic rendering can provide insight
into the structure and significance of the correlation neighborhoods
as we sweep the variables over different vertex degree curves in the
waterfall plot.

The outline of this paper is as follows. In
Sec. \ref{sec:parcorscreening} we formally define the hub screening
problem. In Sec. \ref{sec:Zscorerep} we present the Z-score
representation for the pseudo-inverse of the sample correlation
matrix. In Sec. \ref{sec:theory} we give an overview of the results
pertaining to phase transition thresholds and establish general limit
theorems for the familywise discovery rates and p-values. Section
\ref{sec:statementasythms} gives the formal statements of the results
in the paper. The proofs of these results are given in the appendix.
 In Sec. \ref{sec:experiments} we validate the
theoretical predictions by simulation and illustrate the application
of hub screening to gene microarray data.

\section{Hub screening framework}
\label{sec:parcorscreening}

Let the
$p$-variate $\bX=[X_{1},\ldots, X_{p}]^T$ have mean $\mathbf \mu$ and
non-singular $p\times p$ covariance matrix $\mathbf \Sigma$. We will
often assume that $\bX$ has an elliptically contoured density:
$f_{\bX}(\bx)= g\left((\bx-\mathbf \mu)^T {\mathbf \Sigma}^{-1}
(\bx-\mathbf \mu)\right)$ for some non-negative strictly decreasing
function $g$ on $\Reals^+$.  The correlation matrix and the partial
correlation matrix are defined as $\mathbf \Gamma =
\bD_{\Sigma}^{-1/2} \mathbf \Sigma \bD_{\Sigma}^{-1/2}$ and $\mathbf
\Omega = \bD_{\Sigma^{-1}}^{-1/2} \mathbf \Sigma^{-1}
\bD_{\Sigma^{-1}}^{-1/2} $, respectively, where for a square matrix
$\bA$, $\bD_{A}=\diag(\mathbf A)$ denotes the diagonal matrix obtained
from $\bA$ by zeroing out all entries not on its diagonal.

Available for observation is a $n \times p$ data matrix $\mathbb X$
whose rows are (possibly dependent) replicates of $\bX$:
$$\mathbb X = [\bX_1,\cdots, \bX_p]=[\bX_{(1)}^T, \cdots, \bX_{(n)}^T]^T,$$
where $\bX_i=[X_{1i}, \ldots, X_{ni}]^T$ and $\bX_{(i)}=[X_{i1},
  \ldots, X_{ip}]$ denote the $i$-th column and row, respectively, of
$\mathbb X$.  Define the sample mean of the $i$-th column
$\ol{X}_i=n^{-1} \sum_{j=1}^n X_{ji}$, the vector of sample means
$\ol{\bX}=[\ol{X}_1,\ldots, \ol{X}_p]$, the $p\times p$ sample
covariance matrix $\bS=\frac{1}{n-1} \sum_{i=1}^n
(\bX_{(i)}-\ol{\bX})^T(\bX_{(i)}-\ol{\bX})$, and the $p \times p$
sample correlation matrix $$\bR=\bD_{\bS}^{-1/2}\bS\bD_{\bS}^{-1/2}.$$
For a full rank sample correlation matrix $\bR$ the sample partial
correlation matrix is defined as
$$ \bP = \bD_{R^{-1}}^{-1/2} \bR^{-1} \bD_{R^{-1}}^{-1/2} .$$ In the
case that $\bR$ is not full rank this definition must be
modified. Several methods have been proposed for regularizing the
inverse of a rank deficient covariance including shrinkage and
pseudo-inverse approaches \cite{Schaefer&Strimmer:SAGMB05}. In this
paper we adopt the pseudo-inverse approach and define the sample
partial correlation matrix as
$$ \bP=\bD_{\bR^{\dagger}}^{-1/2}\bR^{\dagger}\bD_{\bR^{\dagger}}^{-1/2},$$
where $\bR^{\dagger}$ denotes the Moore-Penrose pseudo-inverse of $\bR$.

Let the non-negative definite symmetric matrix $\mathbf
\Phi=((\Phi_{ij}))_{i,j=1}^p$ be generic notation for a
correlation-type matrix like $\mathbf\Gamma$, $\mathbf\Omega$,
$\mathbf R$, or $\mathbf P$.  For a threshold $\rho \in [0,1]$ define
${\mathcal G}_\rho({\mathbf \Phi})$ the undirected graph induced by
thresholding $\mathbf \Phi$ as follows.  The graph ${\mathcal
  G}_\rho(\mathbf \Phi)$ has vertex set $\mathcal V = \{1, \ldots,
p\}$ and edge set $\mathcal E = \{e_{ij}\}_{i,j\in \{1, \ldots, p\}:
  i<j}$, where an edge $e_{ij}\in \mathcal E$ exists in ${\mathcal
  G}_\rho({\mathbf \Phi})$ if $|\Phi_{ij}| \geq \rho$. The degree of the
$i$-th vertex of ${\mathcal G}_\rho({\mathbf \Phi})$ is $\left|\{j\neq
i: |\Phi_{ij}| \geq \rho\}\right|$, the number of edges that connect to
$i$.  We call ${\mathcal G}_0({\mathbf \Phi})$ the population
correlation or partial correlation graph \cite{Cox&Wermuth:96} depending on whether $\mathbf \Phi$
is defined as $\mathbf \Gamma$ or $\mathbf \Omega$. Likewise we call
${\mathcal G}_\rho({\mathbf \Phi})$ the empirical correlation or partial correlation graph
depending on whether  $\mathbf \Phi = \bR$ or $\mathbf \Phi=\bP$.
%

A $p\times p$ matrix is said to be row-sparse of degree $k$, called
the {\it sparsity degree}, if no row contains more than $k+1$ non-zero
entries.  When $\mathbf \Phi$ is row-sparse of degree $k$ its graph
${\mathcal G}_\rho({\mathbf \Phi})$ has no vertex of degree greater
than $k$.  A special case is a block-sparse matrix of degree $k$; a
matrix that can be reduced via row-column permutation to block
diagonal form with a single $k\times k$ block.
%

\subsection{\bf The hub screening problem}


Assume that a single
treatment of the $p$ variables yields data matrix $\mX$ of dimension
$n \times p$. A given vertex $i$ is declared a hub screening discovery at
degree level $\delta$ and threshold level $\rho$ if
\be
\left|\{ j: j\neq i, |\Phi_{ij}| \geq \rho \}\right| \geq \delta,
\label{eq:parcorscreen}
\ee where $\mathbf \Phi$ is equal to $\bR$ for correlation hub screening
or is equal to $\bP$ for partial correlation hub screening.  We denote
by $N_{\delta,\rho}\in \{0,\ldots, p\}$ the total number of hub screening
discoveries.

To be practically useful, we need guidelines for selecting the
threshold screening parameters in (\ref{eq:parcorscreen}).
In the sequel we will develop a
large $p$ asymptotic analysis to address the following issues: 1)
phase transitions in the number of discoveries as a function
of these screening parameters; 2) relations between the false
positive rate of the discoveries and these screening parameters.

\subsection{Z-score representation}
\label{sec:Zscorerep}

\def\mU{{\mathbb U}}
\def\mZ{{\mathbb Z}}
\def\mY{{\mathbb Y}}
\def\mX{{\mathbb X}}
\def\mT{{\mathbb T}}

Define the $n\times p$ matrix of Z-scores associated with the data matrix $\mX$
\be \mT=[\bfT_1, \ldots,\bfT_p]=(n-1)^{-1/2} (\bI_n- n^{-1} \mathbf
1 \mathbf 1^T) \mathbb X \bD^{-1/2}_{\bS} ,
\label{eq:ordAp}
\ee where $\bI_n$ is the $n\times n$ identity matrix and $\mathbf 1 =
    [1,\ldots, 1]^T\in \Reals^{n}$.  
This  Z-score matrix is to be distinguished from the $(n-1)\times p$ 
Z-score matrices $\mU$ and $\mY$, denoted collectively by the notation $\mZ$ in the sequel,  
that are derived from the matrix $\mT$.

We exploit the following Z-score representation of the sample
correlation matrix \be \bR=\mT^T \mT,
\label{eq:ordA}
\ee
and defined a set of equivalent but lower dimensional Z-scores
called U-scores. The U-scores lie in the unit sphere $S_{n-2}$ in
$\Reals^{n-1}$ and are obtained by  projecting away the
rowspace components of $\mT$ in the direction of vector $\mathbf 1$.
More specifically, they are constructed as follows.

\def\bHn{{\mathbf H_{2:n}}} Define the orthogonal $n \times n$ matrix
$\mathbf H=[n^{-1/2}\mathbf 1, \mathbf H_{2:n}]$. The matrix $\bHn$
can be obtained by Gramm-Schmidt orthogonalization and satisfies the properties
$${\mathbf 1}^T
\bH =[\sqrt{n},0,\ldots, 0], \;\; \bHn^T\bHn= \bI_{n-1}. $$
The U-score matrix $\mU=[\bU_1,\ldots,\bU_p]$ is obtained from $\mT$ by the following relation
\be
\left[\begin{array}{c} {\mathbf 0}^T \\ \mU\end{array}\right] = \bH^T \mT.
\label{eq:Udef}
\ee
%
%
\begin{lemmas}
Assume that $n<p$. The Moore-Penrose pseudo-inverse of $\bR$ has the representation
\be \bR^{\dagger} = \mU^T[\mU\mU^T]^{-2} \mU.
\label{lemma:MP}
\ee
\end{lemmas}
The proof of the Lemma simply verifies that $\bQ\defined
\mU^T[\mU\mU^T]^{-2} \mU$ satisfies the Moore-Penrose conditions for
$\bQ$ to be the unique pseudo-inverse of $\bR$: 1) the matrices
$\bQ\bR$ and $\bR\bQ$ are symmetric; 2) $\bR\bQ\bR=\bR$; and
3) $\bQ\bR\bQ=\bQ$ \cite{Golub&VanLoan:89}.

Define the $(n-1) \times p $ matrix of partial correlation Z-scores
$\mY=[\bY_1,\ldots,\bY_p]$, $\bY_i \in S_{n-2}$:
\be
\mY= [\mU\mU^T]^{-1} \mU\bD_{\mU^T[\mU\mU^T]^{-2}\mU}^{-1/2}.
\label{eq:parcorscore}\ee
With this definition Lemma \ref{lemma:MP} gives the following Z-score
representation for the sample partial correlation matrix
\be \bP = \mY^T \mY.
\label{eq:Pscore}
\ee

\section{Overview of results}
\label{sec:theory}

%

\subsection{Phase transitions for hub discoveries}
\label{sec:phasetransitions}

There is a phase transition in the number of discovered variables as a
function of the applied screening threshold(s).
This phase transition
threshold, which we call $\rho_{c,\delta}$, is such that if the
screening threshold $\rho$ decreases below $\rho_{c,\delta}$, the mean
number $E[N_{\delta,\rho}]$ of hub discoveries of degree $\delta$  abruptly increases to
$p$. An expression for the critical threshold is obtained from the
asymptotic expression (\ref{eq:ENdef}) for the mean given in
Prop. \ref{prop:parcor1} in the Appendix
\be \rho_{c,\delta}=\sqrt{1-(c_{n,\delta}(p-1))^{-2\delta/(\delta(n-2)-2)}},
\label{eq:rhocrit} \ee
where $c_{n,\delta} = \ra_n \delta J_{p,\delta}$
and $\ra_n=2B((n-2)/2,1/2)$ with
$B(i,j)$ the beta function. The unknown covariance matrix $\mathbf \Sigma$ influences
$\rho_{c,\delta}$ only through the quantity $J_{p,\delta}=J(\ol{f_{\bZ_{\ast_1},
    \ldots, \bZ_{\ast_{\delta+1}}}})$, defined in (\ref{eq:Jdef}),
which is a measure of average
$(\delta+1)$-order dependency among the Z-scores $\{\bZ_i\}_{i=1}^p$.
For large $p$ the constant $c_{n,\delta}$ depends only weakly on $p$ and the critical
threshold increases to $1$ at rate $O((p-1)^{-2\delta/(\delta(n-2)-2)})$, which
is close to logarithmic in $p$ for large $n$ but much faster than logarithmic
for small $n$.

When the rows of $\mX$ are i.i.d. elliptically distributed
 and $\mathbf \Sigma$ is
block-sparse of degree $k$
then, from Prop. \ref{prop:prep}
\be
J(\ol{f_{\bZ_{\ast_1}, \ldots, \bZ_{\ast_{\delta+1}}}})=1+O\left((k/p)^{\gamma_{\delta}}\right),
\label{eq:Jpdasy}
\ee
where $\gamma_{\delta}=\delta+1$ for correlation hub screening and $\gamma_{\delta}=1$
for partial correlation hub screening.
When $p$ is large $O\left((k/p)^{\gamma_{\delta}}\right)$ goes to zero and
 the phase transition threshold $\rho_{c,\delta}$ becomes independent of $\mathbf \Sigma$.

\subsection{Asymptotic p-values on hub discoveries}
\label{sec:pvalues}

As $p\rightarrow \infty$ Prop. \ref{prop:parcor1} in the Appendix states that $N_{\delta,\rho}$ is
dominated by an asymptotically Poisson random  variable $N_{\delta,\rho}^*$ if:
1) $\rho=\rho_p$ increases to one at a prescribed rate depending on $n$;
2) the sparsity degree increases only as $k=o(p)$; and 3) the dependency coefficient
$\|\Delta_{p,n,k,\delta}\|_1$, defined in (\ref{eq:deltapij}), converges to zero.
This guarantees that $P(N_{\delta,\rho}>0)$ converges to the Poisson probability
$1-\exp(-\Lambda_{\delta,\rho})$ where $\Lambda_{\delta,\rho}$
is the rate parameter of $N_{\delta,\rho}^*$. The rate of convergence is provided in
Prop. \ref{prop:parcor} along with a finite $p$ approximation
for $\Lambda_{\delta,\rho}$
$$\Lambda_{\delta,\rho}=\lambda_{\delta,\rho}J(\ol{f_{\bZ_{\ast_1}, \ldots,    \bZ_{\ast_{\delta+1}}}}),$$
with
\be
\lambda_{\delta,\rho}=\lim_{p\rightarrow \infty} p{p-1 \choose \delta}P_0(\rho_p,n)^\delta
,
\label{eq:lambdadef}
\ee
and $P_0(\rho,n)$ is the spherical cap probability defined in (\ref{eq:Podef}).
The parameter
$\lambda_{\delta,\rho}$ does not depend on the underlying distribution of the observations.

Furthermore, when the rows of $\mX$ are i.i.d. elliptically distributed
with covariance that is block-sparse of degree $k$,
Prop. \ref{prop:prep} asserts that (\ref{eq:Jpdasy}) holds and
that the dependency coefficient
$\|\Delta_{p,n,k,\delta}\|_1$ is equal to zero for correlation hub screening
and of order $O(k/p)$ for partial correlation hub screening.
Therefore, when the block-sparse model is posed as the null hypothesis,
Prop. \ref{prop:parcor1} implies that the false positive FWER error rate
can be approximated by
\be P(N_{\delta,\rho}>0)\approx 1-\exp(-\lambda_{\delta,\rho}).
\label{eq:FPR}
\ee

Approximate p-values can also be obtained
under the block-sparse null hypothesis.  Assume that for an
initial threshold $\rho^*\in [0,1]$ the sample correlation or partial correlation
graph ${\mathcal G}_{\rho^*}(\mathbf \Phi)$ has been computed for
$\mathbf \Phi$ equal to $\bR$ or $\bP$.  Consider a vertex $i$ in this
graph that has degree $d_i>0$.
For each value $\delta\in\{1, \ldots, \max_{i=1,\ldots, p} d_i\}$
of the degree threshold $\delta$ denote by $\rho_i(\delta)$
the maximum value of the correlation threshold $\rho$ for which this
vertex maintains degree at least $\delta$ in ${\mathcal G}_{\rho}(\mathbf \Phi)$.
$\rho_i(\delta)\geq \rho^*$ and is equal to the sample correlation between $\bX_i$ and
its $\delta$-th nearest neighbor. We define the approximate p-value
associated with discovery $i$ at degree level $\delta$ as %
\be pv_\delta(i)\approx 1-\exp(-\lambda_{\delta,\rho_i(\delta)}).
\label{eq:PPvalue}
\ee
The quantity (\ref{eq:PPvalue}) is an approximation to
the probability of observing
at least one vertex with degree greater than or equal to
$\delta$ in ${\mathcal G}_\rho({\mathbf \Phi})$ under the nominal block-sparse covariance model.

The accuracy of the approximations of
false positive rate (\ref{eq:FPR}) and p-value
(\ref{eq:PPvalue})  are specified by the bound
(\ref{eq:Pvlim}) given in Prop. \ref{prop:parcor}.
Corollary \ref{cor:parcor}  provides
rates of convergence under the assumptions that
$p(p-1)^{\delta} (1-\rho^2)^{(n-2)/2}=O(1)$ and
the rows of $\mX$ are i.i.d.  with sparse covariance. For example, assume that
the covariance is block-sparse  of degree $k$. If $k$ does not grow with $p$
then the rate of convergence of
$P(N_{\delta,\rho}>0)$ to its Poisson limit is no worse
than
$O(p^{-1/\delta})$ for $\delta>n-3$.  On the other hand, if $k$ grows
with rate at least $O(p^{1-\alpha})$, for
$\alpha=\min\{(\delta+1)^{-1},(n-2)^{-1}\}/\delta$, the rate of
convergence is no worse than $O\left(k/p\right)$. This latter bound
can be replaced by $O\left((k/p)^{\delta+1}\right)$ for correlation
hub screening under the less restrictive assumption that the
covariance is row-sparse.

More generally the combination of
Prop. \ref{prop:parcor} and the assertions (Prop. \ref{prop:prep})
that $J(\ol{f_{\bZ_{\ast_1}, \ldots, \bZ_{\ast_{\delta+1}}}})=1+O\left((k/p)^{\gamma_{\delta}}\right)$ and
$\|\Delta_{p,n,k,\delta}\|_1\leq O(k/p)$ yields
\ben
\left|P(N_{\delta,\rho}>0)-\left(1-\exp(-\lambda_{\delta,\rho})\right)\right|\leq
\left\{\begin{array}{cc} O\left( \max\left\{(k/p)^{\gamma_{\delta}},
p^{-(\delta-1)/\delta}
(k/p)^{\delta-1},p^{-1/\delta},(1-\rho)^{1/2}\right\}\right), &
\delta>1 \\ O\left( \max\left\{(k/p)^{\gamma_{\delta}}, (k/p)^2,
p^{-1},(1-\rho)^{1/2}\right\}\right), & \delta=1 \end{array}\right.
\label{eq:Pvlim_main} .
\een
The terms $(k/p)^{\gamma_{\delta}}$, $p^{-(\delta-1)/\delta}
(k/p)^{\delta-1},p^{-1/\delta}$ and $(1-\rho)^{1/2}$ respectively quantify the
contribution of errors due to: 1)  insufficient sparsity in the covariance or,
equivalently, the correlation graph; 2)
excessive dependency among neighbor variables in this graph; 3)
poor convergence of $E[N_{\delta,\rho}]$; and 4)
inaccurate  mean-value approximation of the integral representation of
$\lim_{p \rightarrow \infty}E[N_{\delta,\rho}]$
by (\ref{eq:thetamoms}). 
One of these
terms will dominate depending on the regime of operation.
For example, specializing to
partial correlation hub screening ($\gamma_{\delta}=\delta+1$), if $\delta>1$
and $O\left(p^{1-(\delta+1)/(2\delta)}\right)\leq k \leq o(p)$ then
$(k/p)^{\delta+1}>p^{-(\delta-1)/\delta} (k/p)^{\delta-1}$ and the
deficiency in the Poisson probability
approximation will not be the determining factor
on convergence rate.
%
%

\subsection{P-value trajectories and waterfall plots}
\label{sec:waterfallplot}

As defined above, the approximate p-values provide a useful
statistical summary.  Rank ordering and thresholding the list
$\{pv_\delta(i)\}_{i=1}^p$ of p-values at any level $\alpha \in [0,1]$
yields the set of vertices of degree at least $\delta$ that would pass
a test of significance at false positive FWER level
$\alpha$. Additional useful information can be gleaned by graphically
rendering the aggregate lists of p-values as described below.

Specifically, assume that correlation screening generates an associated family
of graphs $\left\{ {\mathcal G}_{\rho}(\mathbf \Phi)\right\}_{\rho\in
  [0,1]}$.  Let $d_{\max}=\max_{i=1,\ldots,p} d_i$ be the maximum
discovered degree in the initial graph $ {\mathcal G}_{\rho^*}(\mathbf
\Phi)$.  We define the {\it waterfall plot of p-values} as the family
of curves, plotted against the thresholds $\rho_i(\delta)$,
indexed by $\delta \in \left\{1, \delta, d_{\max}\right\}$
where the $\delta$-th curve is formed from the (linearly interpolated)
ranked ordered list of p-values $\{pv_\delta(i_j)\}_{j=1}^{p}$,
$pv_\delta(i_1) \geq \ldots \geq pv_\delta(i_p)$ (see
Fig. \ref{fig:waterfallplot}).

A useful
visualization of the evolution of vertex neighborhoods over the family
$\left\{ {\mathcal G}_{\rho}(\mathbf \Phi)\right\}_{\rho\in [0,1]}$ is
the {\it p-value trajectory} over the waterfall plot.
This trajectory is defined as the ordered list
$\{pv_\delta(i)\}_{\delta=1}^{d_{max}}$ defined for a given
vertex $i$.  All p-value trajectories start at the outermost curve
(curve associated with $\delta=1$) on the waterfall plot and extinguish
at some inner curve (associated with increasing $\delta>1$).
Vertices with the tightest large neighborhoods will tend to have long
trajectories that start in the middle of the outer curve and
extinguish at a curve deep in the waterfall plot, e.g., the
variable labeled ARRB2 in Fig. \ref{fig:waterfallplot}, while vertices
with the tightest small neighborhoods will tend to have short
trajectories that start  near the extremal point of the outer curve,
e.g., the variable labeled CTAG2 in Fig. \ref{fig:waterfallplot} whose trajectory extinguishes
near the bottom right corner of waterfall plot.

\section{Main theorems}
\label{sec:statementasythms}

The asymptotic theory for hub discovery in correlation and partial graphs is presented in the form of three propositions and one corollary. Prop. \ref{prop:parcor} gives a general bound on
the finite sample approximation error associated with the
approximation of the mean and probability of discoveries given in
Prop. \ref{prop:parcor1}.
The results of Props. \ref{prop:parcor} and \ref{prop:parcor1} apply to general
random matrices of the form $\mZ^T\mZ$ where the $p$ columns of $\mZ$
lie on the unit sphere $S_{n-2}\subset \Reals^{n-1}$ and, in view of
(\ref{eq:Udef}) and (\ref{eq:Pscore}), they provide a unified theory
of hub screening for correlation graphs and partial correlation graphs.
Corollary \ref{cor:parcor}
specializes the bounds
presented in Prop. \ref{prop:parcor} 
to the case of sparse correlation graphs using Prop. \ref{prop:prep}.



For $\delta\geq 1$, $\rho \in [0,1]$,
and $\mathbf \Phi$ equal to the sample
correlation matrix $\bR$ or the sample partial covariance matrix $\bP$
we recall the definition of $N_{\delta,\rho}$ as the number of
vertices of degree at least $\delta$ in ${\mathcal G}_{\rho}(\mathbf
\Phi)$. Define
$\tilde{N}_{\delta,\rho}$ as the count of
the number of groups of $\delta$ mutually coincident edges in
${\mathcal G}_{\rho}(\mathbf \Phi)$
\footnote{$\tilde{N}_{\delta,\rho}$ is equivalent to
  the number of subgraphs in ${\mathcal G}_{\rho}$ that are isomorphic
  to a star graph with $\delta$ vertices.}.
In the sequel we will use the key property that $N_{\delta,\rho}=0$ if and only if
$\tilde{N}_{\delta,\rho}=0$.

For $\delta \geq 1$, $\rho \in [0,1]$, and $n>2$ define
\be
\Lambda =\xi_{p,n,\delta,\rho}J(\ol{f_{\bZ_{\ast_1}, \ldots, \bZ_{\ast_{\delta+1}}}})
\label{eq:Lambdadef}
\ee
where
\be
\xi_{p,n,\delta,\rho}=p{p-1\choose \delta} P_0^{\delta}
\label{eq:xidef},
\ee
$P_0=P_0(\rho,n)$ is defined in (\ref{eq:Podef}), $J$ is given in (\ref{eq:Jdef}),  and $\ol{f_{\bfZ_{\ast_1}, \ldots,
\bfZ_{\ast_{\delta+1}}}}$ is the average joint density given in
(\ref{eq:olfdef2}).

We also define the following quantity needed
for the bounds of Prop. \ref{prop:parcor}
\be
\eta_{p,\delta}=p^{1/\delta}(p-1) P_0.
\label{eq:etadef}
\ee Note that
$\xi_{p,n,\delta,\rho}/\eta_{p,\delta}^{\delta}=(\ra_n/(n-2))^{\delta}/\delta!$
to order $O\left(\max\{p^{-1}, 1-\rho\}\right)$, where
$\ra_n=(2\Gamma((n-1)/2))/(\sqrt{\pi}\Gamma((n-2)/2))$. Let
$\varphi(\delta)$ be the function equal to $1$ for $\delta>1$
and equal to $2$ for $\delta=1$.

\begin{propositions}
Let $\mZ=[\bZ_1, \ldots, \bZ_p]$ be a $(n-1)\times p$ random matrix
with $\bZ_i \in S_{n-2}$. Fix integers $\delta$ and $n$ where
$\delta\geq 1$ and $n>2$. Let the joint density of any subset of the
$\bZ_i$'s be bounded and differentiable. Then, with $\Lambda$ defined in (\ref{eq:Lambdadef}),
\be
\left|E[N_{\delta,\rho}] - \Lambda\right| \leq
O\left( \eta_{p,\delta}^{\delta}
\max\left\{\eta_{p,\delta} p^{-1/\delta} ,(1-\rho)^{1/2}\right\}\right)
\label{eq:ENlim}
\ee
Furthermore, let $N^*_{\delta,\rho}$
be a Poisson distributed random variable with rate
$E[N^*_{\delta,\rho}]=\Lambda/\varphi(\delta)$.
If $(p-1)P_0\leq 1$ then, for any integer $k$, $1\leq k \leq p$,
\be
&&\left|P(N_{\delta,\rho}>0)-P(N^*_{\delta,\rho} >0)\right|\leq \label{eq:Pvlim} \\
&&\hspace{0.5in} \left\{\begin{array}{cc}
O\left(\eta_{p,\delta}^{\delta}
\max\left\{\eta_{p,\delta}^{\delta}
\left(k/p\right)^{\delta+1},
Q_{p,k,\delta}, \|\Delta_{p,n,k,\delta}\|_1,p^{-1/\delta},(1-\rho)^{1/2}\right\}\right), & \delta>1 \\
O\left(\eta_{p,1}
\max\left\{\eta_{p,1}
\left(k/p\right)^{2}, \|\Delta_{p,n,k,1}\|_1, p^{-1},(1-\rho)^{1/2}\right\}\right), & \delta=1
\end{array}\right. ,
\nonumber
\ee
with $Q_{p,k,\delta}=\eta_{p,\delta}^{\delta-1}
p^{-(\delta-1)/\delta}\left(k/p\right)^{\delta-1}$ and
$\|\Delta_{p,n,k,\delta}\|_1$ defined in (\ref{eq:Deltapdefavg}).
\label{prop:parcor}
\end{propositions}
%
%
%

The proof of the above proposition is given in the Appendix.  The
Poisson-type limit (\ref{eq:Poissonconv}) is established by showing
that the count $\tilde{N}_{\rho,\delta}$ of the number of
groups of $\delta$ mutually coincident edges in $\mathcal G_\rho$
converges to a Poisson random variable with rate $\Lambda/\varphi(\delta)$.

\begin{propositions}
Let $\rho_p\in[0,1]$ be a sequence converging to one as $p\rightarrow
\infty$ such that $p^{1/\delta}(p-1)(1-\rho_p^2)^{(n-2)/2} \rightarrow
e_{n,\delta}\in (0,\infty)$. Then
\be
\lim_{p\rightarrow\infty}E[N_{\delta,\rho_p}] = \kappa_{n,\delta} \;  \lim_{p\rightarrow \infty}
J(\ol{f_{\bZ_{\ast_1}, \ldots, \bZ_{\ast_{\delta+1}}}}),
\label{eq:ENdef}
\ee
where $\kappa_{n,\delta}=\left(e_{n,\delta} \ra_n/(n-2)\right)^\delta/\delta!$.
Assume that $k=o(p)$ and that the weak dependency condition
$\lim_{p\rightarrow\infty} \|\Delta_{p,n,k,\delta}\|_1=0$ is satisfied. Then
\be
P(N_{\delta,\rho_p}>0)\rightarrow
1-\exp(-\Lambda/\varphi(\delta)).
\label{eq:Poissonconv}
\ee
\label{prop:parcor1}
\end{propositions}

%
%

The proof of Prop. \ref{prop:parcor1} is an immediate and obvious consequence of
Prop. \ref{prop:parcor} and is omitted.

Propositions \ref{prop:parcor} and \ref{prop:parcor1} are general results
that apply to both correlation hub and partial correlation hub screening
under a wide range of conditions.  Corollary \ref{cor:parcor}
specializes these results to the case of sparse covariance and
i.i.d. rows of $\mX$ having elliptical distribution. These are
standard conditions assumed in the literature on graphical models.


\begin{corollaries}
In addition to the hypotheses of Prop. \ref{prop:parcor1} assume that
$n>3$ and that the rows of $\mathbb X$ are i.i.d. elliptically
distributed with a covariance matrix $\mathbf \Sigma$ that is
row-sparse of degree $k$. Assume that $k$ grows as
$O\left(p^{1-\alpha}\right) \leq k\leq o(p)$ where
$\alpha=\min\left\{(\delta+1)^{-1},(n-2)^{-1}\right\}/\delta $. Then,  for
correlation hub screening 
the
asymptotic approximation error in the limit (\ref{eq:Poissonconv})
is upper bounded by $O\left( (k/p)^{\delta+1} \right)$. 
Under the additional assumption that the covariance is block-sparse, for
partial correlation hub screening  this error
is upper bounded by
$O\left( k/p \right)$.
\label{cor:parcor}
\end{corollaries}


The proof of Corollary \ref{cor:parcor} is given in the Appendix. 
The  proposition below specializes these results to sparse covariance.

\begin{propositions}
Let $\mX$ be a $n\times p$ data matrix whose rows are i.i.d.
 realizations of an elliptically distributed
$p$-dimensional vector $\bX$ with mean parameter $\mathbf \mu$ and
covariance parameter $\mathbf \Sigma$.  Let $\mU=[\bU_1,\ldots,
  \bU_p]$ be the matrix of correlation Z-scores (\ref{eq:Udef}) and $\mY=[\bY_1,
  \ldots, \bY_p]$ be the matrix of partial correlation Z-scores (\ref{eq:parcorscore})
defined in Sec. \ref{sec:Zscorerep}.  Assume that the covariance
matrix $\mathbf \Sigma$ is block-sparse of degree $q$. Then
the pseudo-inverse partial correlation matrix $\bP=
\mY^T\mY$ has the representation
\be \bP=\mU^T\mU(1+O(q/p)) .
\label{Prep}\ee
Let $\bZ_i$ denote $\bU_i$ or $\bY_i$ and assume that for $\delta\geq 1$
the joint density of any distinct set of Z-scores $\bU_{i_1},\ldots, \bU_{i_{\delta+1}}$ is bounded and differentiable over $S_{n-2}^{\delta+1}$.
Then the $(\delta+1)$-fold average function $J$ (\ref{eq:olfdef2}) and the dependency
coefficient $\Delta_{p,n,k,\delta}$ (\ref{eq:Deltapdefavg}) satisfy
\be
J(\ol{f_{\bZ_{\ast_1}, \ldots, \bZ_{\ast_{\delta+1}}}})=1+O\left((q/p)^{\gamma_{\delta}}\right) ,
\label{eq:prep2}
\ee
\be
\|\Delta_{p,n,k,\delta}\|_1 = \left\{\begin{array}{cc} O\left((q/p)\right) ,& \varphi=1 \\
0, & \varphi=0 \end{array} \right.
\label{eq:prep3}
\ee
where $\varphi=0$ and $\varphi=1$ for correlation and partial correlation hub screening, respectively, and
$\gamma_{\delta}=\varphi+(\delta+1)(1-\varphi)$.
\label{prop:prep}
\label{prop:parcorblocksparse}
\end{propositions}

\noindent{\it Proof of Proposition \ref{prop:prep}:}

The proof of Proposition \ref{prop:prep} is given in the Appendix. \qed


\section{Experiments}
\label{sec:experiments}

\subsection{Numerical simulation study}
\label{sec:simulations}

\begin{figure}
\begin{center}
\includegraphics[height=5cm]{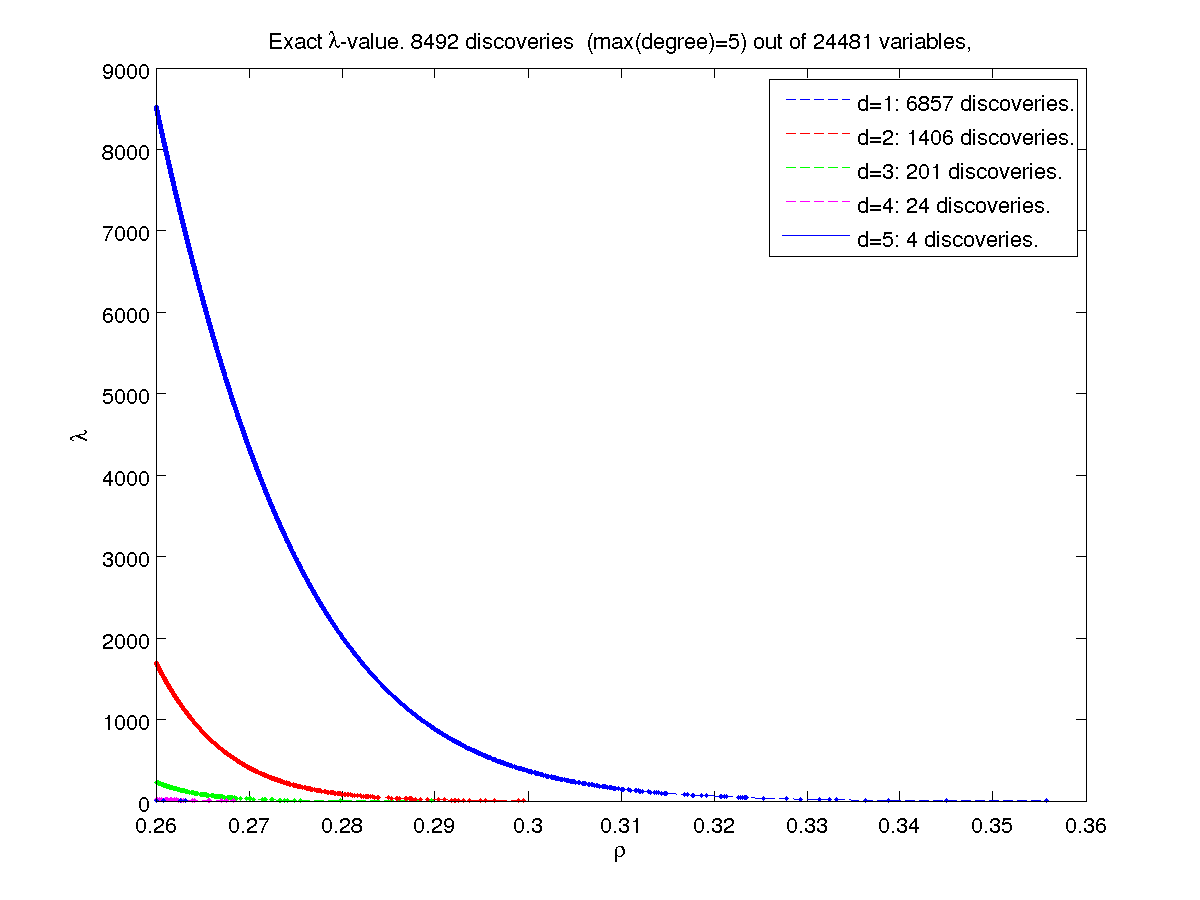} \; \includegraphics[height=5cm]{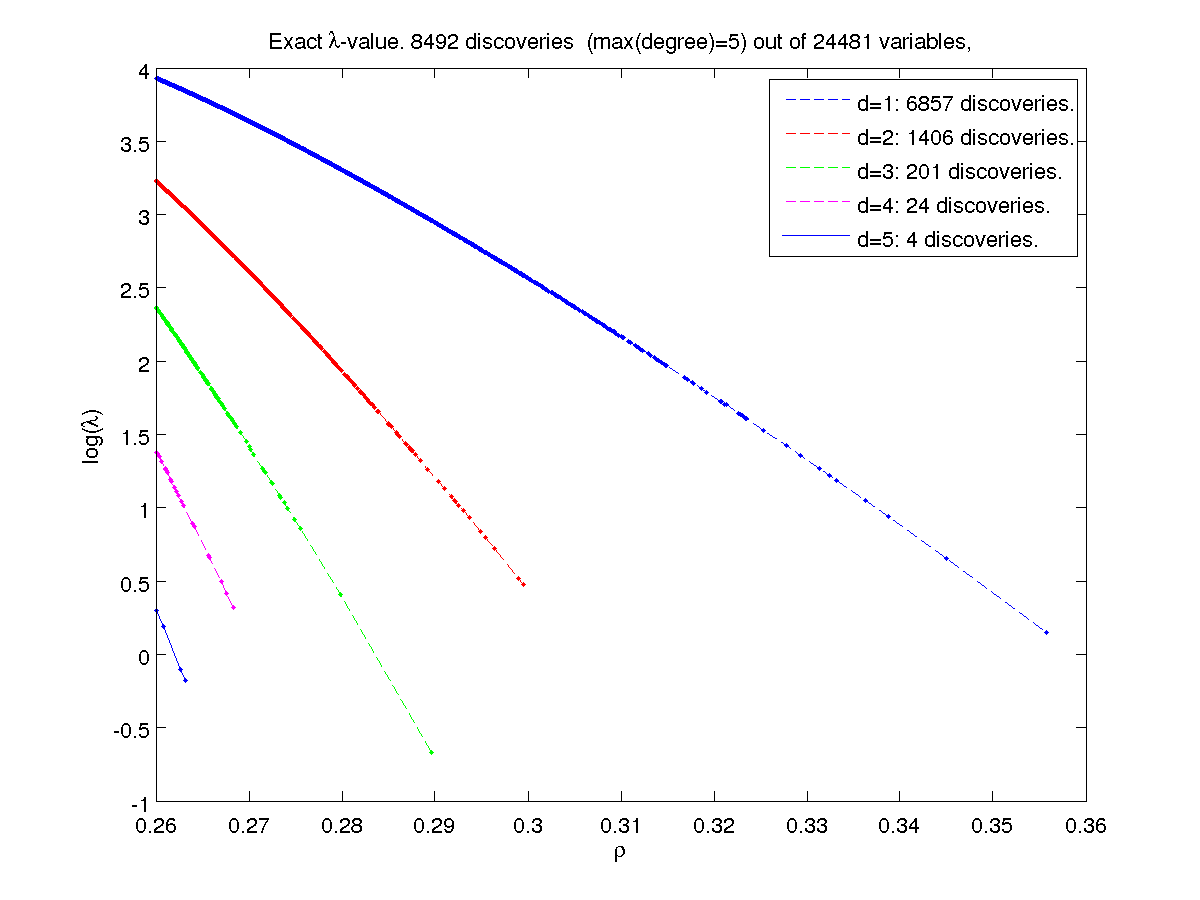}
\end{center}
\caption{Waterfall plots of p-values for partial correlation hub screening of a
  sham version of the NKI dataset \cite{chang2005robustness}.  The
  data matrix $\mX$ has $n=266$ rows and $p=24,481$ columns and is
  populated with i.i.d. zero mean and unit variance Gaussian variables
  (${\mathbf \Sigma}=\bI$). Dots on the curves correspond to
  unnormalized log p-values (called $\lambda$-values) of discovered
  variables whose partial correlations exceeed the initial screening
  threshold $\rho^*=0.26$.  Each curve is indexed over a particular
  vertex degree threshold parameter $\delta$ ranging from $\delta=1$
  to $\delta=5$, the maximum vertex degree found.
  Legends indicate the total number of variables
  discovered having vertex degree $d$.  Left: plot of
  $\lambda_{\delta,\rho_i(\delta)}(i)=-\log(1-pv_\delta(i))$.
  Right: same
  as right panel but $\lambda$-values are plotted on log scale. }
\label{fig:pv1}
\end{figure}

Figure \ref{fig:pv1} shows the waterfall plot of partial correlation hub p-values
 for a sham measurement matrix with i.i.d. normal
entries that emulates the NKI experimental data studied in the next
subsection. There are $n=266$ samples and $p=24,481$ variables in this sham.
Using (\ref{eq:rhocrit}) with parameter $c_{n,\delta} = \ra_n \delta $,
the critical phase transition threshold  on discoveries with positive vertex degree
 was determined to be  $\rho_{1,c}=0.296$.
For purposes of illustration of the fidelity of our theoretical
predictions we used an initial screening threshold equal to
$\rho^*=0.26$. As this is a sham, all discoveries are false positives.

To illustrate the fidelity of the theoretical predictions waterfall
plots of $p$-values (\ref{eq:PPvalue}) are shown in Fig. \ref{fig:pv1}.
For clarity, the figure shows the {\it $\lambda$-value} defined
as $\lambda_{\delta,\rho_i(\delta)}= -\log(1-pv_{\delta}(i)) $.  When
presented in this manner the leftmost point of each curve in the
waterfall plot occurs at approximately $(\rho^*,
E[N_{\delta,\rho^*}])$, as can be seen by comparing the second and
third columns of Table 1. 
This table demonstrates strong
agreement between the predicted (mean) number of partial correlation
hub discoveries and the actual number of discoveries for a single
realization of the data matrix. The realization shown in the table is representative 
of the many simulations performed by us.

\begin{table}
\label{table:sham}
\begin{center}
\begin{tabular}{|c||c|c|}  \hline
observed degree &  \# predicted ($E[N_{\delta,\rho^*}]$) & \# actual ($N_{\delta,\rho^*}$)  \\ \hline
$d_i\geq \delta =1$ & 8531  & 8492 \\ \hline
$d_i\geq \delta =2$ & 1697 & 1635 \\ \hline
$d_i\geq \delta =3$ & 234  & 229 \\ \hline
$d_i\geq \delta =4$ & 24  & 28 \\ \hline
$d_i\geq \delta =5$ & 2 & 4 \\ \hline
\end{tabular}
\caption{Fidelity of the predicted (mean) number of false positives and
  the observed number of false positives in the realization of the sham
NKI dataset experiment shown in Fig. \ref{fig:pv1}}.
\end{center}
\end{table}

\subsection{Parcor screening of NKI dataset}

The Netherlands Cancer Institute (NKI) dataset \cite{chang2005robustness}
contains data from Affymetrix GeneChips collected from 295 subjects
who were diagnosed with early stage breast cancer.
In Peng {\it et al}  \cite{Peng&etal:JASA09} a graphical lasso method
for estimating the partial correlation graph was proposed and was applied to
this dataset.  Here we apply partial correlation hub screening.

\begin{figure}
\begin{center}
\includegraphics[height=10cm]{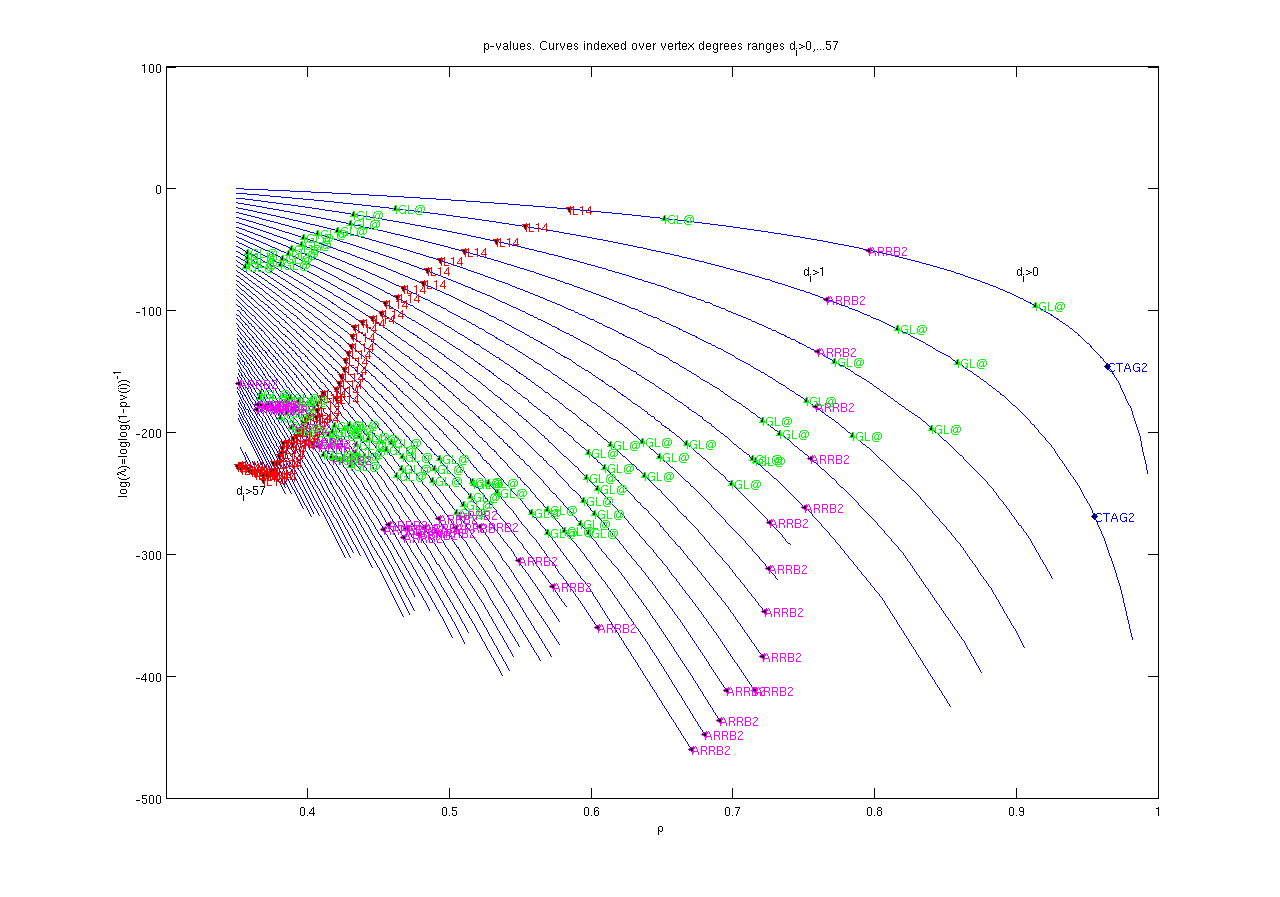}
\end{center}
\caption{Waterfall plot of p-values for NKI gene expression dataset of
  \cite{chang2005robustness} plotted in terms of
  $\log\log(1-pv_\delta(i))^{-1}$.  The genes plotted correspond to
vertices of positive degree in the initial partial correlation graph with threshold $\rho^*=0.35$.
  Each curve indexes the p-values for a particular degree threshold
$\delta$ and a
  gene is on the curve if its degree $d_i$ in the initial graph is
  greater than or equal to $\delta$.  The discovered vertex degree ranges from
  $1$ to $58$ (last dot labeled IL14 at bottom left). The
  p-value trajectories across vertex degree $\delta$ are indicated for several genes
  of interest. Note that all three  Affymetrix array replicates of IGL$@$ were discovered. }
\label{fig:waterfallplot}
\end{figure}

\begin{figure}
\begin{center}
\includegraphics[height=10cm]{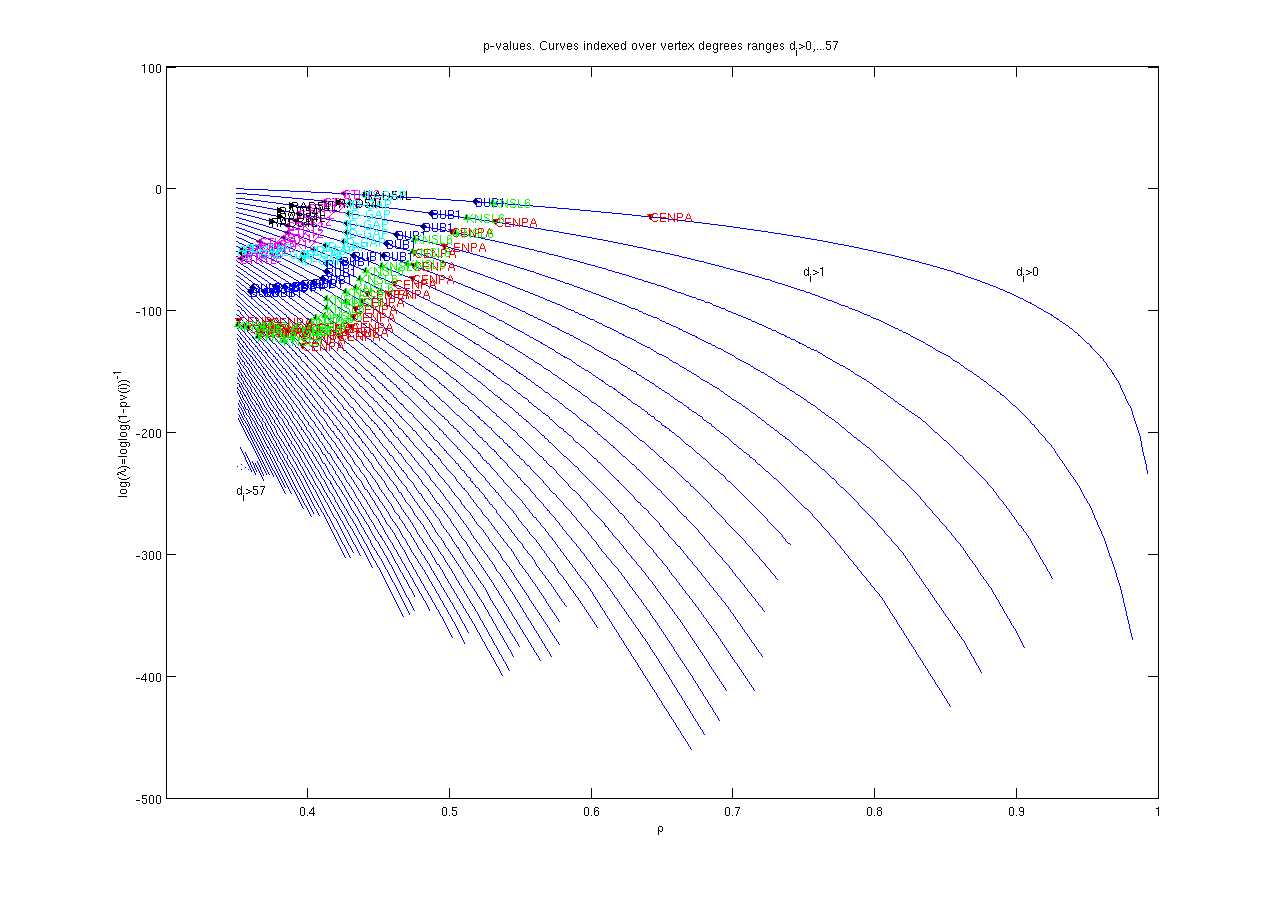}
\end{center}
\caption{Same as Fig. \ref{fig:waterfallplot} but showing the p-value
  trajectories of $6$ of the hub genes ('BUB1' 'CENPA' 'KNSL6' 'STK12'
  'RAD54L' 'ID-GAP') reported in Peng {\it et al}
  \cite{Peng&etal:JASA09}. These are the genes reported in Table 4 of
  \cite{Peng&etal:JASA09} that have unambiguous annotation on the
  GeneChip array.
}
\label{fig:waterfallplot2}
\end{figure}

As in Peng $\etal$ \cite{Peng&etal:JASA09} we only used a subset of
the available GeneChip samples. Specifically, since 29 of the 295
GeneChips had variables with missing values, only 266 of the
them were used in our analysis.  Each GeneChip sample in the NKI
dataset contains the expression levels of 24,481 genes.  Peng
  $\etal$ \cite{Peng&etal:JASA09} applied univariate Cox regression to
reduce the number of variables to 1,217 genes prior to applying their
sparse partial correlation estimation ({\tt space}) method.  In
contrast, we applied our partial correlation hub screening procedure directly to all
24,481 variables. 

An initial threshold $\rho^*=0.35>\rho_{1,c}=0.296$ was selected.
Figure \ref{fig:waterfallplot}
illustrates the waterfall plot of p-values of all discovered variables.
Note in particular the very high level of significance of
certain variables at the lower extremities of the p-value curves.
According to NCBI Entrez several of the most
statistically significant discovered genes on these strands have been
related to breast cancer, lymphoma, and immune response, e.g. ARRB2
({\it Arrestin, Beta 2}), CTAG2 ({\it Cancer/testis antigen}), IL14
({\it Interleukin}), and IGL$@$ ({\it Immunoglobin alpha}). The
p-value trajectories (colored labels) across different values of $\delta$ of these
four genes is illustrated in the figure.  Note that some genes are
statistically significant only at low vertex degree (CTAG2) while
others retain high statistical significance across all vertex degrees
(IGL$@$).  Fig \ref{fig:waterfallplot2} is the same plot with the
trajectories of the $6$ unambiguously annotated hub genes given in
Table 4 of Peng \etal \cite{Peng&etal:JASA09}.  While these $6$
genes do not have nearly as high p-values, or as high partial
correlation, as compared to other genes shown in Fig. \ref{fig:waterfallplot} their
p-values are still very small; less than $10^{-25}$.

\section{Conclusions}

We treated the problem of screening for variables that are strongly
correlated hubs in a correlation or partial correlation graph when $n\ll p$
and $p$ is large.  The proposed hub screening procedure thresholds the
sample correlation or the pseudo-inverse of the sample correlation
matrix using Z-score representations of the
correlation and partial correlation matrices.
For large
$p$ and finite $n$ asymptotic limits that specify the probability of
false hub discoveries were established. These limits were used to
obtain expressions for phase transition thresholds and p-values
under the assumption of a block-sparse covariance matrix.  To
illustrate the wide applicability of our hub screening results we
applied it to the NKI breast cancer gene expression dataset.

\bc
\bf Acknowledgements
\ec

Alfred Hero was partially supported by National Science Foundation
Grant CCF 0830490 and DIGITEO. Bala Rajaratnam was supported in part
by NSF Grants DMS-05-05303, DMS-09-06392 and grant
SUFSC08-SUSHSTF09-SMSCVISG0906.




\begin{thebibliography}{10}

\bibitem{anandkumar2009detection}
A.~Anandkumar, L.~Tong, and A.~Swami.
\newblock {Detection of Gauss--Markov random fields with nearest-neighbor
  dependency}.
\newblock {\em Information Theory, IEEE Transactions on}, 55(2):816--827, 2009.

\bibitem{arratia1990poisson}
R.~Arratia, L.~Goldstein, and L.~Gordon.
\newblock {Poisson approximation and the Chen-Stein method}.
\newblock {\em Statistical Science}, 5(4):403--424, 1990.

\bibitem{Bickel&Levina:AOS08_2}
P.J. Bickel and E.~Levina.
\newblock Covariance regularization via thresholding.
\newblock {\em Annals of Statistics}, 34(6):2577--2604, 2008.

\bibitem{chang2005robustness}
H.Y. Chang, D.S.A. Nuyten, J.B. Sneddon, T.~Hastie, R.~Tibshirani,
  T.~S{\o}rlie, H.~Dai, Y.D. He, L.J. Van't~Veer, H.~Bartelink, et~al.
\newblock {Robustness, scalability, and integration of a wound-response gene
  expression signature in predicting breast cancer survival}.
\newblock {\em Proceedings of the National Academy of Sciences}, 102(10):3738,
  2005.

\bibitem{Cox&Wermuth:96}
D.R. Cox and N.~Wermuth.
\newblock {\em {Multivariate dependencies: models, analysis and
  interpretation}}.
\newblock Chapman \& Hall/CRC, 1996.

\bibitem{Dempster:Biometrics72}
A.P. Dempster.
\newblock {Covariance selection}.
\newblock {\em Biometrics}, 28(1):157--175, 1972.

\bibitem{friedman2007sic}
J.~Friedman, T.~Hastie, and R.~Tibshirani.
\newblock {Sparse inverse covariance estimation with the graphical lasso}.
\newblock {\em Biostatistics}, 2007.

\bibitem{friedman2008sparse}
J.~Friedman, T.~Hastie, and R.~Tibshirani.
\newblock {Sparse inverse covariance estimation with the graphical lasso}.
\newblock {\em Biostatistics}, 9(3):432, 2008.

\bibitem{gill2010statistical}
R.~Gill, S.~Datta, and S.~Datta.
\newblock {A statistical framework for differential network analysis from
  microarray data}.
\newblock {\em BMC bioinformatics}, 11(1):95, 2010.

\bibitem{GoldsteinSmith:74}
M.~Goldstein and A.~F.~M. Smith.
\newblock Ridge-type estimators for regression analysis.
\newblock {\em Journal of the Royal Statistical Society. Series B
  (Methodological)}, 36(2):pp. 284--291, 1974.

\bibitem{Golub&VanLoan:89}
G.~H. Golub and C.~F. {Van Loan}.
\newblock {\em Matrix Computations (2nd Edition)}.
\newblock The Johns Hopkins University Press, Baltimore, 1989.

\bibitem{Hero&Rajaratnam:Paper1_Techreport}
A.O. Hero and B.~Rajaratnam.
\newblock Large scale correlation screening.
\newblock {\em Technical Report, Department of Statistics and EECS, University
  of Michigan \& Department of Statistics, Stanford University.
  http://arxiv.org/abs/1102.1204}, 2010.

\bibitem{Jegou2010product}
H.~J{\'e}gou, M.~Douze, and C.~Schmid.
\newblock {Product quantization for nearest neighbor search}.
\newblock {\em IEEE Transactions on Pattern Analysis and Machine Intelligence},
  2010.

\bibitem{Kramer&etal:BMC09}
Nicole Kr\"amer, Juliane Sch\"afer, and Anne-Laure Boulesteix.
\newblock Regularized estimation of large-scale gene association networks using
  graphical gaussian models.
\newblock {\em BMC Bioinformatics}, 10(384):1--24, 2009.

\bibitem{liu2010feedback}
Y.~Liu, V.~Chandrasekaran, A.~Anandkumar, and A.S. Willsky.
\newblock {Feedback message passing for inference in gaussian graphical
  models}.
\newblock In {\em Information Theory Proceedings (ISIT), 2010 IEEE
  International Symposium on}, pages 1683--1687. IEEE, 2010.

\bibitem{Peng&etal:JASA09}
J.~Peng, P.~Wang, N.~Zhou, and J.~Zhu.
\newblock {Partial correlation estimation by joint sparse regression models}.
\newblock {\em Journal of the American Statistical Association},
  104(486):735--746, 2009.

\bibitem{pihur2008reconstruction}
V.~Pihur, S.~Datta, and S.~Datta.
\newblock {Reconstruction of genetic association networks from microarray data:
  a partial least squares approach}.
\newblock {\em Bioinformatics}, 24(4):561, 2008.

\bibitem{rajaratnam2008flexible}
B.~Rajaratnam, H.~Massam, and C.M. Carvalho.
\newblock {Flexible covariance estimation in graphical Gaussian models}.
\newblock {\em Annals of Statistics}, 36:2818--2849, 2008.

\bibitem{Rothman&etal:08}
A.~J. Rothman, P.J. Bickel, E.~Levina, and J.~Zhu.
\newblock Sparse permutation invariant covariance estimation.
\newblock {\em Electronic Journal of Statistics}, 2:494--515, 2008.

\bibitem{Schaefer&Strimmer:SAGMB05}
J.~Schaefer and K~Strimmer.
\newblock A shrinkage approach to large-scale covariance matrix estimation and
  implications for functional genomics.
\newblock {\em Statist. Appl. Genet. Mol. Biol.}, 4(32), 2005.

\bibitem{wiesel2010covariance}
A.~Wiesel, Y.C. Eldar, and A.O. Hero.
\newblock {Covariance estimation in decomposable Gaussian graphical models}.
\newblock {\em Signal Processing, IEEE Transactions on}, 58(3):1482--1492,
  2010.

\end{thebibliography}

\newpage
\section{Appendix}

This appendix contains two subsections. Section \ref{sec:notation} gives the necessary definitions. Section \ref{sec:asythms} gives proofs of the theory given in Sec. \ref{sec:statementasythms} .

\subsection{Notation, Preliminaries and Definitions}
\label{sec:notation}


\begin{itemize}
\item
$\mX$: $n \times p$ matrix of observations.
\item
$\mZ=[\bZ_1, \ldots , \bZ_p]$: $(n-1)\times p$ matrix of
correlation (or partial correlation) $Z$-scores $\{\bZ_i\}_i$
associated with $\mX$.
\item
$\mZ^T\mZ$: $p\times p$ sample correlation matrix $\bR$ (if $\mZ=\mU$) or sample partial correlation matrix $\bP$ (if $\mZ=\mY)$  associated with $\mX$.
\item
$\bZ_i^T \bZ_j$: sample  correlation (or partial correlation) coefficient,
the $i,j$-th element of $\mZ^T\mZ$.
\item
$\rho\in [0,1]$: screening threshold applied to  matrix $\mZ^T\mZ$.
\item
$r=\sqrt{2(1-\rho)}$: spherical cap radius parameter.
\item
$S_{n-2}$: unit sphere in $\Reals^{n-1}$.
\item
$a_n=|S_{n-2}|$: surface area of $S_{n-2}$.
\item
${\mathcal G}_0({\mathbf \Phi})$:  graph associated with population correlation
matrix  $\mathbf \Phi=\mathbf \Gamma$ or partial correlation matrix $\mathbf
\Phi=\mathbf \Omega$. An edge in ${\mathcal G}_0({\mathbf \Phi})$ corresponds to
 a non-zero entry of $\mathbf \Phi$.
\item
${\mathcal G}_\rho={\mathcal G}_\rho({\mathbf \Phi})$:  graph
 associated with thresholded sample correlation matrix  $\mathbf \Phi=\bR$ or partial correlation matrix $\mathbf \Phi=\bP$.
Specifically, the edges
of ${\mathcal G}_\rho({\mathbf \Phi})$ are specified by the non-diagonal
entries of $\mZ^T\mZ$ whose magnitudes exceed level $\rho$.
\item
$d_i$: observed degree of vertex $i$ in ${\mathcal G}_\rho({\mathbf \Phi})$, $\mathbf \Phi\in \{\bR,\bP\}$.
\item
$\delta$:  screening threshold for vertex degrees
in ${\mathcal G}_\rho({\mathbf \Phi})$, $\mathbf \Phi\in \{\bR,\bP\}$.
\item
$k$: upper bound on vertex degrees of ${\mathcal G}_0({\mathbf \Phi})$,
$\mathbf \Phi\in \{\mathbf \Gamma,\mathbf \Omega\}$.
\item
$N_{\delta,\rho}$: generic notation for  the number of correlation hub discoveries ($N_{\delta,\rho}(\bR)$) or partial correlation hub discoveries ($N_{\delta,\rho}(\bP)$) of degree $d_i\geq \delta$ in  ${\mathcal G}_\rho(\bR)$, or ${\mathcal G}_\rho(\bP)$, respectively.
\item
$\tilde{N}_{\delta,\rho}$ counts the number of subsets of $\delta$
  mutually coincident edges in 
  ${\mathcal G}_{\rho}$.
\item
$A(r,\bz)$: the union of two anti-polar spherical cap regions in $
  S_{n-2}$ of radii $r =\sqrt{2(1-\rho)}$ centered at points $-\bz$
  and $\bz$.
\item
$P_0$: probability that a uniformly distributed vector $\bU\in S_{n-2}$ falls in $A(r,\bz)$
\be
\begin{split}
P_0 = \; P_0(\rho,n)=& \; \ra_n \int_\rho^1\left(1-u^2\right)^{\frac{n-4}{2} } du \\
   =& \; (n-2)^{-1} \ra_n (1-\rho^2)^{(n-2)/2}(1+O(1-\rho^2)) ,
\end{split}
\label{eq:Podef}
\ee
where 
$\ra_n=2B((n-2)/2,1/2)$ and $B(l,m)$ is the beta function.
\end{itemize}

For given integer $k$,
$0\leq k <p$, and $\Phi$ either the population correlation matrix
$\mathbf \Gamma$ or the population partial correlation matrix $\mathbf
\Omega$ define
\be \calN_k(i)=\argmax_{j_1 \neq \cdots \neq
  j_{\min(k,\ol{d}_i)}} \sum_{l=1}^{\min(k,\ol{d}_i)} |\Phi_{ij_l}|
\label{eq:Nki},
 \ee where $\ol{d}_i$ denotes the degree of vertex $i$ in ${\mathcal
   G}_0({\mathbf \Phi})$ and the maximization is over the range of  distinct
$j_l\in \{1, \ldots, p\}$ that are not equal to $i$.
When $k\geq \ol{d}_i$ these are
 the indices of the $\ol{d}_i$ neighbors of vertex $i$ in ${\mathcal
   G}_0({\mathbf \Phi})$. When $k<\ol{d}_i$ these are the subset of
the $k$-nearest neighbors ($k$-NN) of vertex $i$.  For the sequel it will be
 convenient to define the following vector valued indexing variable:
 $\vec{i}=(i_0, \ldots, i_{\delta})$, where $0<\delta\leq p$ and $i_0, \ldots, i_{\delta}$
 are distinct integers in $\{1, \ldots, p\}$.  With this index denote by
 $\bZ_{\vec{i}}$  the set of $\delta+1$ Z-scores
 $\{\bZ_{i_j}\}_{j=0}^\delta$.

Define the set of complementary $k$-NN's of $\bZ_{\vec{i}}$ as
$\bZ_{A_k(\vec{i})}=\{\bZ_l: l \in A_k(\vec{i})\}$, where \be A_{k}(\vec{i}) =
\left(\cup_{l=0}^\delta\calN_k(i_l)\right)^c -\{\vec{i}\}
\label{eq:Akij},
\ee
with $A^c$ denoting set complement of set $A$. The complementary
$k$-NN's include vertices outside of the
$k$-nearest-neighbor regions of the set of points $\bZ_{\vec{i}}$.

Define the $\delta$-fold leave-one-out average of the density, a function of $i$,
$f_{\bZ_{i_1}, \ldots, \bZ_{i_\delta}, \bZ_{i}}$:
\be &&\ol{f_{\bZ_{\ast_1-i}, \ldots, \bZ_{\ast_\delta-i}, \bZ_{i}}}
(\bfz_1, \ldots, \bfz_{\delta}, \bfz_i) \label{eq:avgfubivpersistent}
\label{eq:olfdef1}\\
&&\hspace{0.5in}=2^{-d} \sum_{s_1, \ldots ,s_\delta \in \{-1,1\}} {p-1 \choose \delta}^{-1}
\sum_{i_1 \neq \cdots \neq i_\delta \neq i}^{p}
f_{\bZ_{i_1},\ldots, \bZ_{i_\delta},\bZ_{i}}(s_1 \bfz_1, \ldots, s_\delta\bfz_{\delta}, \bfz_i),
 \nonumber
\ee
where in the inner summation, indices $i_1,\ldots, i_\delta$ range over $\{1, \ldots, p\}$.
Also define the $(\delta+1)$-fold average of the same density
\be
&&\ol{f_{\bZ_{\ast_1}, \ldots,  \bZ_{\ast_{\delta+1}}}}
(\bfz_1, \ldots, \bfz_{\delta}, \bfz_i) \label{eq:olfdef2} \\ &&\hspace{0.5in}=
p^{-1}\sum_{i=1}^{p} \left(\half \ol{f_{\bZ_{\ast_1-i}, \ldots, \bZ_{\ast_\delta-i}, \bZ_{i}}}
(\bfz_1, \ldots, \bfz_{\delta}, \bfz_i) +
\half\ol{f_{\bZ_{\ast_1-i}, \ldots, \bZ_{\ast_\delta-i}, \bZ_{i}}}
(\bfz_1, \ldots, \bfz_{\delta}, -\bfz_i)\right).
 \nonumber
\ee

For any data matrix $\mZ$ define the
dependency coefficient between the columns $\bZ_{\vec{i}}$ and
their complementary $k$-NN's
\be
\Delta_{p,n,k,\delta}(\vec{i})= \left \|(f_{\bZ_{\vec{i}}|\bZ_{A_k(\vec{i})}}-f_{\bZ_{\vec{i}}})/f_{\bZ_{\vec{i}}}
\right\|_{\infty},
\label{eq:deltapij}
\ee
and the average of these  coefficients is
\be
\|\Delta_{p,n,k,\delta}\|_1= \left(p{p-1 \choose \delta}\right)^{-1}
\sum_{i_0=1}^p \sum_{i_1< \ldots < i_\delta}
\Delta_{p,n,k,\delta}(\vec{i}) .
\label{eq:Deltapdefavg}
\ee%
where the second sum is indexed over $i_1,\ldots, i_\delta \neq i_0$.

The coefficient (\ref{eq:Deltapdefavg})
quantifies weak  dependence of the Z-scores.  If, for all $i$,
$\bZ_{\vec{i}}$ and its complementary $k$-NN neighborhood variables are
independent then $\|\Delta_{p,n,k,\delta}\|_1=0$.  When the rows of $\mX$ are
i.i.d. and elliptically distributed, and $\mZ=\mU$ are the standard
correlation $Z$-scores, then a sufficient condition for
$\|\Delta_{p,n,k,\delta}\|_1=0$ is that ${\mathcal G}_0({\mathbf \Phi})$
have no vertex of
degree greater than $k$ or, equivalently, that the dispersion matrix
$\mathbf \Sigma$ be row sparse of degree $k$.

Finally, for arbitrary
joint density $f_{\bZ_1,\ldots,\bZ_\delta}(\bz_1,\ldots,\bz_\delta)$ on
$S_{n-2}^\delta={\large \times}_{i=1}^\delta S_{n-2}$,  define
\be J(f_{\bZ_1,\ldots,\bZ_\delta}) = |S_{n-2}|^{\delta-1} \int_{S_{n-2}} f_{\bZ_1,\ldots,\bZ_\delta}(\bz, \ldots,\bz) d\bz .
\label{eq:Jdef}
\ee

\subsection{Proofs of  theorems}
\label{sec:asythms}

\noindent{\it Proof of Prop. \ref{prop:parcor}: }

With $\phi_i=I(d_i\geq \delta)$
we have  $N_{\delta,\rho} =
\sum_{i=1}^p \phi_i$.  Define $\phi_{ij}=I\left(\bfZ_j\in
A(r,\bfZ_i)\right)$ the indicator of the presence of an edge in
${\mathcal G}_\rho({\mathbf \Phi})$ between vertices $i$ and $j$,
where $A(r,\bfZ_i)$ is the union
of two antipolar caps  in $S_{n-2}$ of radius
$r=\sqrt{2(1-\rho)}$ centered at $\bfZ_i$ and $-\bfZ_i$,
respectively. Then $\phi_i$ and $\phi_{ij}$ have the explicit relation
\be
\phi_i&=&\sum_{l=\delta}^{p-1} \sum_{\vec{k}\in\breve{\calC}_i(p-1,l)} \prod_{j=1}^l \phi_{ik_j}
\prod_{m=l+1}^{p-1} (1-\phi_{i k_m})
\label{eq:phirep}
\ee
where we have defined the index vector $\vec{k}=(k_1, \ldots, k_{p-1})$ and the set
$$\breve{\calC}_i(p-1,l)=\{\vec{k}: k_1< \ldots< k_l, k_{l+1} < \ldots < k_{p-1}
\; k_j \in \{1, \ldots, p\}-\{i\}, k_j\neq k_{j'}\}.$$
The inner summation in (\ref{eq:phirep}) simply sums over the set of distinct indices not equal to $i$ that index all
${p-1 \choose l}$ different types of products $\prod_{j=1}^l \phi_{ik_j} \prod_{m=l+1}^{p-1} (1-\phi_{i k_m})$.
Subtracting $
\sum_{\vec{k}\in\breve{\calC}_i(p-1,\delta)}  \prod_{j=1}^\delta \phi_{ik_j}$
from both sides of (\ref{eq:phirep})
\be
&&\phi_i-\sum_{\vec{k} \in \breve{\calC}_i(p-1,\delta)}
\prod_{j=1}^\delta \phi_{ik_j}\\ &&\hspace{0.2in}=
\sum_{l=\delta+1}^{p-1} \sum_{\vec{k}\in\breve{\calC}_i(p-1,l)}  \prod_{j=1}^l \phi_{ik_j}\prod_{m=l+1}^{p-1} (1-\phi_{i k_m})
\\
&&\hspace{0.3in}+\sum_{\vec{k}\in\breve{\calC}_i(p-1,l)}
\sum_{m=l+1}^{p-1}(-1)^{m-l} \sum_{k_{l+1} < \ldots <k_m} \prod_{j=1}^l \phi_{ik_j}  \prod_{n=l+1}^m \phi_{ik_n}
\label{eq:phirep2}
\ee
where, in the last line we have used the  expansion
$$\prod_{m=l+1}^{p-1} (1-\phi_{i k_m})=
1+ \sum_{m=l+1}^{p-1}(-1)^{m-l} \sum_{k_{l+1} < \ldots <k_{m}} \prod_{n=l+1}^m \phi_{ik_n}.$$

The following simple asymptotic representation will be useful in the sequel. For any $i_1, \ldots, i_k \in \{1, \ldots , p\}$,
$i_1 \neq \cdots \neq i_k \neq i$, $k\in \{1, \ldots, p-1\}$,
\be E\left[\prod_{j=1}^k \phi_{i i_j} \right] &=&\int_{S_{n-2}} d\bfv \int_{A(r,\bfv)} d\bfu_{1} \cdots
\int_{A(r,\bfv)} d\bfu_{k}\;
f_{\bfU_{i_1},\ldots, \bfU_{i_k},\bfU_i}(\bfu_1,
\cdots,\bfu_k,\bfv)
\label{eq:thetamoms}\\
&\leq& P_0^k a_n^k M_{k|1}
\label{eq:thetamomsineq}
\label{eq:old714}
\ee
with $P_0=P_0(\rho,n)$ defined in (\ref{eq:Podef}), $a_n=|S_{n-2}|$, and
\be
M_{k|1} &=&\max_{i_1 \neq \cdots \neq i_{k+1}}\left\| f_{\bfZ_{i_1},\ldots, \bfZ_{i_k}|\bfZ_{i_{k+1}}}\right\|_{\infty}
\label{eq:Mk1def},
\ee  
The following simple generalization of (\ref{eq:old714})
to arbitrary product indices $\phi_{ij}$ will also be needed
\be
E\left[\prod_{l=1}^m \phi_{i_lj_l}\right] \leq  P_0^m a_n^m M_{|Q|},
\label{eq:thetamoms2}
\ee
where $Q=$unique$(\{i_l,j_l\}_{l=1}^m)$ is the set of unique indices
among the distinct pairs $\{(i_l,j_l)\}_{l=1}^m$ and $M_{|Q|}$ is a bound on the joint density of $\bZ_{Q}$.

Define the random variable
\be\theta_i = {p-1 \choose \delta}^{-1}
\sum_{\vec{k} \in \breve{\calC}_i(p-1,\delta)} \prod_{j=1}^\delta \phi_{ik_j}.
\label{eq:thetadefnew}
\ee
We show below that for sufficiently large $p$
\be
\left|E[\phi_i]-{p-1 \choose \delta}
E[\theta_i]\right| &\leq& \gamma_{p,\delta} ((p-1) P_0)^{\delta+1},
\label{eq:proof1}
\ee
where $\gamma_{p,\delta}=\max_{\delta+1\leq l<p} \{a_n^l M_{l|1}\}
\left(e-\sum_{l=0}^\delta \frac{1}{l!}\right) \left(1+(\delta!)^{-1}\right)$ and
$M_{l|1}$ is a least upper bound on any $l$-dimensional joint
density of the variables $\{\bZ_i\}_{j\neq i}^p$ conditioned on $\bZ_i$.

To show inequality (\ref{eq:proof1}) take expectations of (\ref{eq:phirep2})
and apply the bound (\ref{eq:old714}) 
to obtain
\be
&&\left|E[\phi_i]-{p-1 \choose \delta} E[\theta_i]\right| \nonumber \\
&&\hspace{0.3in} \leq
\left|\sum_{l=\delta+1}^{p-1} {p-1 \choose l}  P_0^l a^lM_{l|1}
+{p-1 \choose \delta}  \sum_{l=1}^{p-1-\delta} {p-1-\delta \choose l} P_0^{\delta+l} a_n^{\delta+l}
M_{\delta+l|1} \right| \nonumber \\
&&\hspace{0.3in} \leq A(1+(\delta!)^{-1}),
\label{eq:proof2}
\ee
where
$$A=\sum_{l=\delta+1}^{p-1} {p-1 \choose l} P_0^l a^lM_{l|1}.$$
The line (\ref{eq:proof2})
follows from the identity
${p-1-\delta \choose l}{p-1 \choose \delta} ={p-1 \choose l+\delta} (\delta!)^{-1}$
and a change of index in the second summation on the previous line.
Since $(p-1)P_0<1$
\ben |A|
&\leq& \max_{\delta+1\leq l<p} \{a_n^l M_{l|1}\} \sum_{l=\delta+1}^{p-1} {p-1 \choose l}
((p-1)P_0)^l \\
&\leq& \max_{\delta+1\leq l<p} \{a_n^l M_{l|1}\}\left(e-\sum_{l=0}^\delta \frac{1}{l!}\right)
((p-1)P_0)^{\delta+1}.
\een

Application of the mean-value-theorem to the integral representation (\ref{eq:thetamoms}) yields
\be
\left|E[\theta_i]- P_0^\delta J(\ol{f_{\bZ_{\ast_1-i}, \ldots, \bZ_{\ast_\delta-i},\bZ_i}})
\right| &\leq& \tilde{\gamma}_{p,\delta} ((p-1)P_0)^\delta r,
\label{eq:proof3}
 \ee
where 
$\tilde{\gamma}_{p,\delta}= 2a_n^{\delta+1} \dot{M}_{\delta+1|1}/\delta!$ and
$\dot{M}_{\delta+1|1}$ is a bound on the norm of the gradient
$$\nabla_{\bz_{i_1}, \ldots, \bz_{i_\delta}}\ol{f_{\bZ_{\ast_1-i}, \ldots, \bZ_{\ast_\delta-i}|\bZ_i}(\bz_{i_1}, \ldots, \bz_{i_\delta}|\bz_i)}.$$
Combining (\ref{eq:proof1})-(\ref{eq:proof3})
and the relation $r=O((1-\rho)^{1/2})$,
\be
\left|E[\phi_{i}]- {p-1 \choose \delta} P_0^{\delta} J(\ol{f_{\bZ_{\ast_1}, \ldots,
            \bZ_{\ast_{\delta+1}}}})\right| &\leq& O\left(((p-1)P_0)^{\delta} \max\left\{ (p-1)P_0, (1-\rho)^{1/2}\right\} \right) .
\label{eq:proof3pp}
\ee
Summing over $i$ and recalling  the definitions (\ref{eq:xidef})  and  (\ref{eq:etadef})
of $\xi_{p,n,\delta,\rho}$ and $\eta_{p,\delta}$,
\be
\left|E[N_{\delta,\rho}]-\xi_{p,n,\delta,\rho}
J(\ol{f_{\bZ_{\ast_1}, \ldots,
            \bZ_{\ast_{\delta+1}}}})\right| &\leq&
O\left(p((p-1)P_0)^{\delta} \max\left\{ (p-1)P_0, (1-\rho)^{1/2}\right\} \right)
\nonumber
\\
&=&O\left(\eta_{p,\delta}^{\delta}
\max\left\{\eta_{p,\delta}
p^{-1/\delta}, (1-\rho)^{1/2}\right\} \right).
\label{eq:proof3p}
\ee
This establishes the bound (\ref{eq:ENlim}).

For the bound (\ref{eq:Pvlim}) we
use the Chen-Stein method \cite{arratia1990poisson}.
The part of the bound (\ref{eq:Pvlim}) that holds for $\delta=1$
was derived in the course of proof
of Prop. 1 in \cite{Hero&Rajaratnam:Paper1_Techreport}. Below we
treat the case $\delta>1$.
Recalling the definition  $\tilde{N}_{\delta,\rho}$ as the number
of subsets of $\delta$ mutually coincident edges
in ${\mathcal G}_{\rho}$, we have the representation:
\be
\tilde{N}_{\delta,\rho}= \sum_{i_0=1}^p \sum_{i_1<\ldots < i_{\delta}}\prod_{j=1}^{\delta}\phi_{i_0i_j}
\label{eq:rep1},
\ee
where the second sum is indexed over $i_1, \ldots, i_\delta \neq i_0$.
For $\vec{i}\defined (i_0,i_1, \ldots, i_{\delta})$ define the
index set $\Beta_{\vec{i}}=\Beta_{i_0,i_1,\ldots, i_{\delta}}=\{(j_0,j_1,\ldots,
j_{\delta}): j_l\in \mathcal N_k(i_l)\cup\{i_l\}, l=0, \ldots, \delta\}\cap
{\mathcal C}^{<}$ where ${\mathcal C}^{<}=\{ (j_0, \ldots, j_\delta): j_0\in \{1,\ldots,p\}, $ $1\leq j_1< \cdots < j_\delta \leq p, j_1, \ldots, j_\delta\neq j_0\}$.  These index the distinct sets of points
$\bZ_{\vec{i}}=\{\bZ_{i_0},\bZ_{i_1}, \ldots, \bZ_{i_{\delta}}\}$ and their respective
$k$-NN's.  Note that $|\Beta_{\vec{i}}|\leq k^{\delta+1}$. 
Identifying
$\tilde{N}_{\delta,\rho}= \sum_{\vec{i}\in {\mathcal C}^{<}}
\prod_{l=1}^{\delta}\phi_{i_0i_l}$ and $N_{\delta,\rho}^*$ a Poisson distributed random
variable with rate $E[\tilde{N}_{\delta,\rho}]$, the Chen-Stein bound
\cite[Thm. 1]{arratia1990poisson} is
\be 2\max_A |P(\tilde{N}_{\delta,\rho} \in A)-P(N_{\delta,\rho}^*\in
A)| \leq b_1+b_2+b_3,
\label{eq:chenstein}
\ee
where
$$b_1 = \sum_{\vec{i}\in {\mathcal C}^{<}}\sum_{\vec{j}\in \Beta_{\vec{i}}}
E\left[\prod_{l=1}^{\delta}\phi_{i_0i_l}\right]E\left[\prod_{m=1}^{\delta}\phi_{j_0j_m}\right],
$$
$$b_2 =
\sum_{\vec{i}\in {\mathcal C}^{<}}\sum_{\vec{j}\in \Beta_{\vec{i}}-\{\vec{i}\}}
E\left[\prod_{l=1}^{\delta}\phi_{i_0i_l}\prod_{m=1}^{\delta}\phi_{j_0j_m}\right],$$
and, for $p_{\vec{i}}=E[\prod_{l=1}^{\delta}\phi_{i_0i_l}]$,
$$
b_3 = \sum_{\vec{i}\in {\mathcal C}^{<}}
E\left[ E\left[\left.\prod_{l=1}^{\delta}\phi_{i_0i_l}-p_{\vec{i}}\right|\phi_{\vec{j}}: \vec{j} \not\in \Beta_{\vec{i}}\right]\right].$$

Over the range of indices in the sum of $b_1$
$E[\prod_{l=1}^{\delta}\phi_{ii_l}]$ is of order $O(P_0^{\delta})$, by (\ref{eq:thetamoms2}),
and therefore
$$b_1\leq O\left(p^{\delta+1} k^{\delta+1} P_0^{2\delta} \right)=
O\left(\eta_{p,\delta}^{2\delta}(k/p)^{\delta+1}\right),$$
which follows from definition (\ref{eq:etadef}).
More
care is needed to bound $b_2$ due to the symmetry relation
$\phi_{ij}=\phi_{ji}$. If in the summations defining $b_2$,
$i_0=j_m$ and $j_0=i_l$ occur for some $l,m$
then there will be a match and $\phi_{i_0i_l}\phi_{j_0j_m}=\phi_{i_0i_l}$. In such case
the summand of $b_2$ will be of lower order than
$O(P_0^{2\delta})$. For example, for the case that $l,m=1$ a match implies
$\phi_{i_0i_1}=\phi_{j_0j_1}$ and, from (\ref{eq:thetamoms2}),
$$E[\prod_{l=1}^{\delta}\phi_{i_0i_l}\prod_{m=1}^{\delta}\phi_{j_0j_m}]=E[\prod_{l=1}^{\delta}\phi_{j_1
    i_l}\prod_{m=2}^{\delta}\phi_{i_1 j_m}] = O\left(P_0^{2\delta-1}\right).$$
Over ${\mathcal C}^{<}$ and $\Beta_{\vec{i}}-\{i\}$
there can be no more than a single match in $b_2$'s summand. For a given match there are at most
$p^{\delta+1}k^{\delta-1}$ summands of reduced order. We conclude that
\ben
b_2&\leq& O\left(p^{\delta+1} k^{\delta+1} P_0^{2\delta} \right)+O\left(p^{\delta+1} k^{\delta-1} P_0^{2\delta-1} \right) \\
&=&
O\left(\eta_{p,\delta}^{2\delta}(k/p)^{\delta+1}\right)+
O\left(\eta_{p,\delta}^{2\delta-1}(k/p)^{\delta-1} p^{-(\delta-1)/\delta} \right) ,
\een
which follows from the relation $p^{2\delta}P_0^{2\delta-1}=(p^{\delta+1}P_0^{\delta})^{2-1/\delta}/p^{(\delta-1)/\delta}$.


Next we bound the term $b_3$ in (\ref{eq:chenstein}). The set
$A_k(\vec{i})=  \Beta_{\vec{i}}^c - \{\vec{i}\}$ indexes
the complementary $k$-NN's of $\bZ_{\vec{i}}$ so that, using the representation   (\ref{eq:thetamoms2}),
\ben
b_3 &=& \sum_{\vec{i}\in {\mathcal C}^{<}}
E\left[E\left[ \left.\prod_{l=1}^{\delta}\phi_{i_0i_l}-p_{\vec{i}}\right|\bZ_{A_k(\vec{i})}\right]\right]
\\
&=&  \sum_{\vec{i}\in {\mathcal C}^{<}}
\int_{S_{n-2}^{|A_k(\vec{i})|}} d\bfz_{A_k(\vec{i})} \left(\prod_{l=1}^{\delta} \int_{S_{n-2}}
d\bfz_{i_0} \int_{A(r,\bfz_{i_0})} d\bfz_{i_l}\right) 
\left(
\frac{f_{\bfZ_{\vec{i}} |\bfZ_{A_k}}(\bfz_{\vec{i}}|\bfz_{A_k(\vec{i})})-
f_{\bfZ_{\vec{i}}}(\bfz_{\vec{i}})}{f_{\bfZ_{\vec{i}}}(\bfz_{\vec{i}})} \right)
f_{\bfZ_{\vec{i}}}(\bfz_{\vec{i}})f_{\bfZ_{A_k(\vec{i})}}(\bfz_{A_k(\vec{i})})
\\
&\leq& O\left(p^{\delta+1}  P_0^{\delta} \|\Delta_{p,n,k,\delta}\|_1\right) = O\left(\eta_{p,\delta}^{\delta}\|\Delta_{p,n,k,\delta}\|_1\right).\een
Observe that,  with $\Lambda=E[N_{\delta,\rho}]$ 
\be \left|P(N_{\delta,\rho} >0)-\left(1-\exp(-\Lambda)\right)\right| &\leq&
\left|P(\tilde{N}_{\delta,\rho} >0)-P(N_{\delta,\rho} >0)\right| \nonumber\\ &&+
\left|P(\tilde{N}_{\delta,\rho} >0)-\left(1-\exp(-E[\tilde{N}_{\delta,\rho}])\right)\right|\nonumber \\
&&+
\left|\exp(-E[\tilde{N}_{\delta,\rho}])-\exp(-\Lambda)\right| \nonumber \\
&\leq & b_1+b_2+b_3 + O\left(\left|E[\tilde{N}_{\delta,\rho}]-\Lambda\right|\right).
\label{eq:allterms}
\ee
Combining the above inequalities on $b_1$, $b_2$ and $b_3$
yields the first three terms in the argument
of the ``max'' on the right side of (\ref{eq:Pvlim}).

It remains to bound the  term $|E[\tilde{N}_{\delta,\rho}]-\Lambda|$.
Application of the mean value theorem to the multiple
integral (\ref{eq:thetamoms2}) gives
\be
\left|E\left[\prod_{l=1}^{\delta}\phi_{ii_l}\right] - P_0^{\delta}
J\left(f_{\bfZ_{i_1},\ldots, \bfZ_{i_{\delta}},\bfZ_{i}}\right)\right| &\leq& O\left(P_0^{\delta}
r \right).
\label{eq:phi_ijdbnd}
\ee
Applying relation (\ref{eq:rep1}) 
%
\be
\left|E[\tilde{N}_{\delta,\rho}] - p{p-1 \choose \delta} P_0^{\delta}
J\left(\ol{f_{\bZ_{\ast_1}, \ldots, \bZ_{\ast_{\delta+1}}}}\right)\right| &\leq& O\left(p^{\delta+1} P_0^{\delta} r\right)
= O\left(\eta_{p,\delta}^{\delta} r\right).
\label{eq:phi_ijdbnd2}
\ee
Combine this with (\ref{eq:allterms}) to obtain the bound (\ref{eq:Pvlim}). This completes the proof of Prop. \ref{prop:parcor}.
\qed

\noindent{\it Proof of Cor. \ref{cor:parcor}: }


For correlation hub screening ($\mZ=\mU$) $\|\Delta_{p,n,k,\delta}\|_1=0$ so it suffices to consider
the other arguments of ``max'' in  the bound (\ref{eq:Pvlim}).
As in the proof of Prop \ref{prop:parcor1},
$(1-\rho_p)^{1/2}= O\left(p^{-(\delta+1)/((n-2)\delta)}\right)$ and we
can merge the last two terms in (\ref{eq:Pvlim}) into the single term
$p^{-\alpha(\delta+1)}=\max\left\{p^{-1/\delta}, (1-\rho_p)^{1/2}\right\} $,
with $\alpha$ defined in the Corollary statement.   Finally, note that $\eta_{p,\delta}=O(1)$ and
$(k/p)^{\delta+1}\geq Q_{p,k,\delta}=(k/p)^{\delta-1}p^{-(\delta-1)/\delta}$
when $k/p\geq p^{-(\delta-1)/(2\delta)}$. Therefore, as $\alpha\leq (\delta-1)/(2\delta)$ when  $n>3$,
we conclude that if $k/p\geq p^{-\alpha}$ all arguments of ``max'' in the bound (\ref{eq:Pvlim}) are dominated by
$(k/p)^{\delta+1}$.
Turning to partial correlation hub screening ($\mZ=\mY$), under the block-sparse covariance assumption
Prop. \ref{prop:prep}  asserts that $\|\Delta_{p,n,k,\delta}\|_1=O(k/p)$  which dominates $(k/p)^{\delta+1}$.
This completes the proof of Cor. \ref{cor:parcor}. \qed

\noindent{\it Proof of Proposition \ref{prop:prep}:}

By block-sparsity, the matrix  $\mU$ of Z-scores can be
partitioned as $\mU
=[\widetilde{\mU},\ol{\mU}]$, where
$\widetilde{\mU}=[\widetilde{\mU}_1,\ldots, \widetilde{\mU}_{q}]$ and
$\ol{\mU}=[\ol{\mU}_1,\ldots, \ol{\mU}_{p-q}]$ are the dependent and
independent columns of $\mU$, respectively. Since the columns of
$\ol{\mU}$'s are i.i.d. and uniform over
the unit sphere $S_{n-2}$, as $p\rightarrow \infty$ we have
\def\olU{\ol{\mU}}
\def\whU{\widetilde{\mU}}
$$(p-q)^{-1} \olU \ \olU^T \rightarrow E[\olU_i \olU_i^T]= (n-1)^{-1} \bI_{n-1}.$$
Furthermore, as the entries of the matrix $q^{-1}\whU\whU^T$ are bounded by 1,
$$ p^{-1} \whU\whU^T = \bO(q/p),$$
where $\bO(u)$ is an $(n-1)\times (n-1)$ matrix
whose entries are of order $O(u)$.
Hence, as $\mU\mU^T = \olU\olU^T + \whU\whU^T$, the pseudo-inverse of $\bR$ has the asymptotic large $p$ representation
%
\be \bR^{\dagger} =\left( \frac{n-1}{p} \right)^2
\mU^T[\bI_{n-1}+\bO(q/p)]^{-2} \mU =\left( \frac{n-1}{p} \right)^2
\mU^T \mU(1+O(q/p)),
\label{Rdaggerrep}
\ee
which establishes (\ref{Prep}).

Define the partition  ${\mathcal C} = \calQ\cup \calQ^c$ of the index set $\calC=\{(i_0,\ldots, i_\delta): 1\leq i_0\neq \cdot \neq i_{\delta}\leq p\}$  where $\calQ=\{(i_0,
\ldots, i_{\delta}): 1\leq i_l \leq q, 1\leq l\leq \delta\} $
is the set of $(d+1)$-tuples restricted to
the dependent columns $\widetilde{\mU}$ of $\mU$.
The  summation representations (\ref{eq:olfdef2}) and
(\ref{eq:Deltapdefavg}) of $J$ and $\|\Delta_{p,n,k,\delta}\|_1$
yield
\be
J(\ol{f_{\bZ_{\ast_1}, \ldots, \bZ_{\ast_{\delta+1}}}})=|\mathcal C|^{-1}
\left(\sum_{\vec{i} \in \mathcal Q}+\sum_{\vec{i} \not \in \mathcal Q}\right)
J(f_{\bZ_{i_0}, \ldots, \bZ_{i_{\delta}}}),\label{eq:Jprep}
\ee
and
\be
\|\Delta_{p,n,k,\delta}\|_1=|\mathcal C|^{-1}
\left(\sum_{\vec{i} \in \mathcal Q}+\sum_{\vec{i} \not \in \mathcal Q}\right)
\Delta_{p,n,k,\delta}(\vec{i}).
\label{eq:Dprep}
\ee
For correlation hub screening ($\mZ=\mU$) $\Delta_{p,n,k,\delta}(\vec{i})=0$ for all $\vec{i}\in \calC$
while, as the set $\{\bU_{i_0},\ldots,\bU_{i_{\delta}}\} $'s
are i.i.d. uniform for $\vec{i}\in \calQ^c$, $J(f_{\bZ_{i_0}, \ldots, \bZ_{i_{\delta}}})=1$ for $\vec{i}\in \calQ^c$. As
$J(f_{\bZ_{i_0}, \ldots, \bZ_{i_{\delta}}})$ is bounded and  $|\calQ|/|\calC|=O\left((q/p)^{\delta+1}\right)$
the relations (\ref{eq:prep2}) and (\ref{eq:prep3}) are
established for the case of correlation screening.

For partial correlation hub screening ($\mZ=\mY$) then, as $\mY=[I_{n-1}+\bO(q/p)]^{-1} \mU$,
the joint densities of
$\mY$ and $\mU$ are related by $f_{\mY}=(1+O(q/p)) f_{\mU}$.
Therefore, over the range $\vec{i} \not \in  \mathcal Q$, the $J$ and
$\Delta_{p,n,k,\delta}$ summands in (\ref{eq:Jprep}) and (\ref{eq:Dprep}) are
of order $1+O(q/p)$ and $O(q/p)$, respectively, which establishes (\ref{eq:prep2}) and (\ref{eq:prep3}) for partial correlation screening.
This completes the proof of Prop. \ref{prop:prep}. \qed


\end{document}